\documentclass{article}%
\usepackage{amsmath}
\usepackage{amsfonts}
\usepackage{amssymb}
\usepackage{graphicx}%
\newtheorem{theorem}{Theorem}
\newtheorem{acknowledgement}[theorem]{Acknowledgement}

\newtheorem{claim}[theorem]{Claim}

\newtheorem{corollary}[theorem]{Corollary}
\newtheorem{criterion}[theorem]{Criterion}
\newtheorem{definition}[theorem]{Definition}

\newtheorem{lemma}[theorem]{Lemma}

\newtheorem{proposition}[theorem]{Proposition}

\newenvironment{proof}[1][Proof]{\noindent\textbf{#1.} }{\ \rule{0.5em}{0.5em}}
\begin{document}

\title{Geometric triangulations and discrete Laplacians on manifolds}
\author{David Glickenstein\\University of Arizona\\and\\Massachusetts Institute of Technology}
\maketitle
\tableofcontents

\section{Introduction\label{introductions}}

In this paper we shall explore Euclidean structures on manifolds which lead to
Laplace operators. Euclidean structures can be introduced on a triangulation
of a manifold by giving each simplex the geometric structure of a Euclidean
simplex. This structure gives the manifold a length space structure in the
same way a Riemannian metric gives a manifold a length structure: the length
between two points is the infimum of the lengths of paths between the two
points. The length of a path is determined by the fact that each simplex it
passes through has the structure of Euclidean space.

The purpose of this paper is to be able to do analysis on the piecewise
Euclidean space. The Laplace operator $\triangle$ is well defined on many
geometric spaces, and is especially important as a natural operator on a
Riemannian manifold and as a generator of Brownian motion. In this paper, we
define a general Euclidean structure called a duality triangulation which not
only allows one to measure length between points and volume of simplices, but
also allows one to describe a geometric dual cell decomposition and the volume
of dual cells. This allows one to define a Laplace operator in a natural way,
which has been applied to fields such as image processing \cite{MDSB}
\cite{Hir} and physics \cite{Mer}.

The duality triangulation structure is very similar to other Euclidean
structures used in both pure and applied math; specifically, we address the
connection to weighted triangulations and Thurston triangulations. In
addition, positivity of volumes of certain duals correspond to Delaunay or
regular triangulations, which are used in a very wide range of applications
from biology to physics to computer graphics.

This paper is organized as follows. We begin in Section
\ref{euclidean structures} with an introduction to Euclidean structures by
recalling the definitions of weighted and Thurston triangulations, introducing
dual triangulations, and relating the three types of triangulations. In
Section \ref{regular triangulations} we discuss regular triangulations and
Delaunay triangulations and consider flip algorithms for constructing regular
and Delaunay triangulations. In Section \ref{laplace} we introduce the Laplace
operator $\triangle$ associated to a given duality triangulation and derive
some of its properties. Finally, in Section \ref{Riemannian} we briefly
discuss the status of piecewise linear Riemannian geometry.

The major new results in this paper are the result on the equivalence of
weighted, Thurston, and duality triangulations in Section
\ref{equivalence of triangulations}, the analysis of flip algorithms in
Section \ref{regular triangulations}, the generalization of Rippa's theorem to
regular triangulations in Section \ref{rippa}, and the definiteness results in
Section \ref{laplace and heat equations}.

Many of the results in this paper were motivated as generalizations of those
described in \cite{BS}.

\section{Euclidean structures\label{euclidean structures}}

\subsection{Basic definitions\label{Euclidean structures basic}}

In this section we shall introduce three types of Euclidean structures:
weighted triangulations, Thurston triangulations, and duality triangulations.
All structures begin with a topological triangulation $\mathcal{T=}\left\{
\mathcal{T}_{0},\mathcal{T}_{1},\ldots,\mathcal{T}_{n}\right\}  $ of an
$n$-dimensional manifold (we shall usually use $n$ to denote the dimension of
the complex in this paper). The triangulation consists of lists of simplices
$\sigma^{k}$, where the superscript denotes the dimension of the simplex, and
$\mathcal{T}_{k}$ is a list of all $k$-dimensional simplices $\sigma
^{k}=\left\{  i_{0},\ldots,i_{k}\right\}  $. We shall often refer to
$0$-dimensional simplices as vertices, $1$-dimensional simplices as edges,
$2$-dimensional simplices as faces or triangles, and $3$-dimensional simplices
as tetrahedra. We shall often denote vertices as $j$ instead of $\left\{
j\right\}  .$ Let $\mathcal{T}_{1}^{+}$ denote the directed edges, where we
distinguish $\left(  i,j\right)  $ from $\left(  j,i\right)  $. When the order
does not matter, we use $\left\{  i,j\right\}  $ to denote an edge. A
triangulation is said to be an $n$-dimensional manifold if a neighborhood of
every vertex is homeomorphic to a ball in $\mathbb{R}^{n}.$ A two-dimensional
manifold is often referred to as a surface. Throughout this paper we will be
dealing exclusively with triangulations of manifolds or parts of manifolds.

In order to give the topological triangulation a geometric structure, each
edge $\left\{  i,j\right\}  $ is assigned a length $\ell_{ij}$ such that for
each simplex in the triangulation there exists a Euclidean simplex with those
edge lengths. We call such an assignment a \emph{Euclidean triangulation}
$\left(  \mathcal{T},\mathcal{\ell}\right)  $, where we think of $\ell$ as a
function
\[
\ell:\mathcal{T}_{1}\rightarrow(0,\infty).
\]
The conditions on $\ell$ include the triangle inequality, but there are
further restrictions in higher dimensions which ensure that the simplices can
be realized as (non-degenerate) Euclidean simplices. The restrictions can be
expressed in terms of the square of volume, which can be expressed as a
polynomial in the squares of the edge lengths by the Cayley-Menger determinant
formula. Each pair of simplices $\sigma_{1}^{n}$ and $\sigma_{2}^{n}$
connected at a common boundary simplex $\sigma^{n-1}$ is called a
\emph{hinge}. In a Euclidean triangulation every hinge can be embedded
isometrically in $\mathbb{R}^{n}$.

Euclidean triangulations have the structure of a distance space with an
intrinsically defined distance. Given any curve $\gamma$ whose length can be
computed on each Euclidean simplex, we can compute the total length of the
curve $L\left(  \gamma\right)  $ as $L\left(  \gamma\right)  =\sum_{\sigma
}L_{\sigma}\left(  \gamma\cap\sigma\right)  $ where $L_{\sigma}\left(
\gamma\cap\sigma\right)  $ is the length of the curve in the simplex $\sigma$
(if the curve intersects the simplex many times, we simply add the
contributions of each piece of the intersection). In particular, we can
consider curves which are differentiable when restricted to each simplex
(these are called piecewise differentiable curves). The intrinsic distance is
defined as
\begin{equation}
d\left(  P,Q\right)  =\inf\left\{  L\left(  \gamma\right)  :\gamma\text{ is a
path from }P\text{ to }Q\right\}  . \label{intrinsic distance}%
\end{equation}
The class of paths can be either taken to be piecewise differentiable or
piecewise linear since length is minimized on piecewise linear paths, as
explained in \cite[Section 2]{Sto}. A path which locally minimizes length is
called a \emph{geodesic} and one which globally minimizes is called a
\emph{minimizing geodesic}.

We are now ready to introduce more structures on Euclidean triangulations.

\subsection{Weighted triangulations\label{weighted triangulations}}

We begin with weighted triangulations.

\begin{definition}
A \emph{weighted triangulation }is a Euclidean triangulation $\left(
\mathcal{T},\ell\right)  $ together with weights
\[
w:\mathcal{T}_{0}\rightarrow\mathbb{R}.
\]

\end{definition}

\begin{figure}
[ptb]
\begin{center}
\includegraphics[
natheight=6.260800in,
natwidth=8.490400in,
height=2.4713in,
width=3.3449in
]
{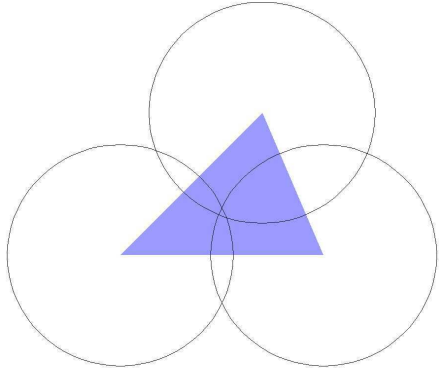}
\caption{A weighted or Thurston triangulation with corresponding circles at
the vertices.}
\label{circletriangle}
\end{center}
\end{figure}

\begin{figure}
[ptb]
\begin{center}
\includegraphics[
natheight=7.776400in,
natwidth=8.490400in,
height=2.9307in,
width=3.1982in
]
{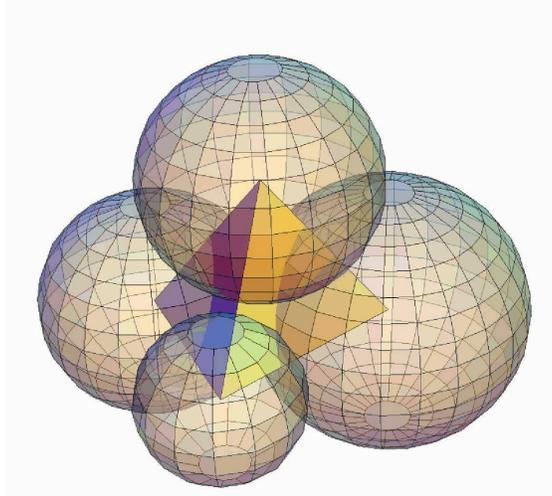}
\caption{A weighted or Thurston triangulation with corresponding spheres at
the vertices.}
\label{spheretetrahedron}
\end{center}
\end{figure}

We think of the weight $w_{i}$ as the square of the radius of a circle
centered at the vertex $i.$ These weighted triangulations are used in the
literature on regular triangulations such as \cite{ES} and \cite{AK}. Thinking
of the weights in this way, in each $n$-dimensional simplex there exists an
$\left(  n-1\right)  $-dimensional sphere which is orthogonal to each of the
spheres centered at the vertices (this means they are perpendicular if they
intersect, or else orthogonal in the sense described in \cite[Section 40]%
{Ped}). In this way, each simplex $\sigma$ has a corresponding center
$C\left(  \sigma\right)  ,$ which is the center of this sphere, and the center
has a weight $w_{C\left(  \sigma\right)  }$ which is the square of the radius
of this sphere. See Figures \ref{circletriangle} and \ref{spheretetrahedron}.%

An important particular case of weighted triangulations is that when $w_{i}=0$
for all vertices $i.$ This is the basis for Delaunay triangulations, but may
not satisfy the Delaunay condition. We shall revisit this in Section
\ref{regular triangulations}.

\subsection{Thurston triangulations\label{thurston triangulations}}

\begin{definition}
A \emph{Thurston triangulation} is a collection $\left(  \mathcal{T}%
,w,c\right)  ,$ where
\begin{align*}
w  &  :\mathcal{T}_{0}\rightarrow\mathbb{R},\\
c  &  :\mathcal{T}_{1}\rightarrow\mathbb{R},
\end{align*}
where $c_{ij}<w_{i}+w_{j}$ and such that the induced lengths
\[
\ell_{ij}=\sqrt{w_{i}+w_{j}-c_{ij}}%
\]
make $\left(  \mathcal{T},\ell\right)  $ into a Euclidean triangulation.
\end{definition}

For a Thurston triangulation, one considers the weight $w_{i}$ to be the
square of the radius $r_{i}$ of a sphere centered at vertex $i,$ just as for
weighted triangulations, and one considers $c_{ij}=2r_{i}r_{j}\cos\left(
\pi-\theta_{ij}\right)  $ where $\theta_{ij}$ is the angle between the spheres
centered at vertices $i$ and $j.$ In this case, one derives the formula for
$\ell_{ij}$ by the law of cosines. By considering $c_{ij}$ instead of
$\theta_{ij},$ we have included some cases where the spheres do not intersect.
These structures were studied by W. Thurston in the context of proving
Andreev's theorem (see \cite{Thurs} and \cite{MR}).

An important special case is that when $c_{ij}=-2r_{i}r_{j}$ (i.e.
$\theta_{ij}=0$). This is the case of a sphere packing on each simplex, since
it corresponds to the spheres being mutually tangent (as in \cite{CR}
\cite{G1} \cite{G2}).

\subsection{Duality triangulations\label{duality triangulations}}

\begin{definition}
A \emph{duality triangulation} is a collection $\left(  \mathcal{T},d\right)
,$ where%
\[
d:\mathcal{T}_{1}^{+}\rightarrow\mathbb{R}%
\]
which satisfies
\begin{equation}
d_{ij}^{2}+d_{jk}^{2}+d_{ki}^{2}=d_{ji}^{2}+d_{ik}^{2}+d_{kj}^{2}
\label{compatibility for local duality}%
\end{equation}
for each $\left\{  i,j,k\right\}  \in\mathcal{T}_{2}$ and such that the
induced lengths
\[
\ell_{ij}=d_{ij}+d_{ji}%
\]
make $\left(  \mathcal{T},\ell\right)  $ into a Euclidean triangulation.
\end{definition}

We think of the weight $d_{ij}$ as representing the portion of the length
$\ell_{ij}$ of edge $\left\{  i,j\right\}  $ which has been assigned to vertex
$i$ while $d_{ji}$ is the portion assigned to vertex $j.$ We thus call them
\emph{local lengths}. The total length of $\left\{  i,j\right\}  $ is the sum
of the contributions $d_{ij}$ from vertex $i$ and $d_{ji}$ from vertex $j.$
Hence each edge is assigned a center $C\left(  \left\{  i,j\right\}  \right)
$ which is distance $d_{ij}$ from vertex $i$ and distance $d_{ji}$ from vertex
$j.$ The condition (\ref{compatibility for local duality}) ensures that for
each triangle $\left\{  i,j,k\right\}  ,$ the perpendiculars to the three
edges through the edge centers meet at one point, which can be called the
center of the triangle, $C\left(  \left\{  i,j,k\right\}  \right)  .$ We shall
soon see that this condition on $2$-dimensional simplices allows us to define
a center for every simplex in the triangulation.

There are two canonical examples which automatically satisfy the condition
(\ref{compatibility for local duality}). One is the case where $d_{ij}$
depends only on $i$ for all edges $\left(  i,j\right)  $ (that is,
$d_{ij}=d_{ik},$ etc.). We call this a \emph{circle} or \emph{sphere packing}
as in \cite{G1}, and the dual comes from the inscripted circle, that is, the
center $C\left(  \left\{  i,j,k\right\}  \right)  $ is the center of the
circle inscribed in $\left\{  i,j,k\right\}  $ in 2D and the center $C\left(
\left\{  i,j,k,\ell\right\}  \right)  $ is the center of the sphere tangent to
each of the edges of the tetrahedron $\left\{  i,j,k,\ell\right\}  $ in 3D.
Another important case is where $d_{ij}=d_{ji}.$ This corresponds to the
center $C\left(  \left\{  i,j,k\right\}  \right)  $ coming from the circle
circumscribed about the triangle $\left\{  i,j,k\right\}  $ and similar for
all higher dimensions.

The structure is called a duality triangulation because the existence of a
center $C\left(  \sigma\right)  $ for each $\sigma$ puts a piecewise-Euclidean
length structure on the dual of the triangulation in such a way that dual
simplices are orthogonal to ordinary simplices. For example, in two
dimensions, if an edge $\left\{  i,j\right\}  $ is part of the two simplices
$\left\{  i,j,k\right\}  $ and $\left\{  i,j,\ell\right\}  ,$ then we can
define the length of the dual edge $\bigstar\left\{  i,j\right\}  $ to be
equal to the distance from the center $C\left(  \left\{  i,j,k\right\}
\right)  $ of the triangle $\left\{  i,j,k\right\}  $ to the center $C\left(
\left\{  i,j\right\}  \right)  $ of the edge $\left\{  i,j\right\}  $ plus the
distance from $C\left(  \left\{  i,j,\ell\right\}  \right)  $ to $C\left(
\left\{  i,j\right\}  \right)  .$ When the hinge is isometrically embedded in
$\mathbb{R}^{2},$ we see that $\bigstar\left\{  i,j\right\}  $ is a straight
line which is perpendicular to the edge $\left\{  i,j\right\}  .$ We shall now
show that this can be done in all dimensions, and no additional restrictions
must be made besides (\ref{compatibility for local duality}) for each triangle.

\begin{proposition}
\label{existence of duals}A duality triangulation in any dimension has unique
centers $C\left(  \sigma^{m}\right)  $ for each simplex $\sigma^{m}$ such that
$C\left(  \sigma^{m}\right)  $ is at the intersection of the $\left(
m-1\right)  $-dimensional hyperplanes through $C\left(  \left\{  i,j\right\}
\right)  $ and perpendicular to $\left\{  i,j\right\}  $ for each $\left\{
i,j\right\}  $ in $\sigma^{m}.$
\end{proposition}

\begin{proof}
We construct the centers $C\left(  \sigma^{m}\right)  $ inductively for
$m$-dimensional simplices. Each pair of $m$-dimensional simplices meeting at
an $\left(  m-1\right)  $-dimensional simplex (a \textquotedblleft
hinge\textquotedblright) can be embedded in $\mathbb{R}^{m}$ as two adjacent
Euclidean simplices. To make the notation more readable, we shall not
distinguish between the embedding of the hinge in $\mathbb{R}^{m}$ and the
hinge as abstract simplices in the piecewise Euclidean manifold. A simplex
$\sigma^{m}$ is assumed to be Euclidean with the assigned edge lengths given
by $\ell_{ij}.$ We now inductively construct the centers of each simplex.
First, $C\left(  \left\{  i\right\}  \right)  =i$ and $C\left(  \left\{
i,j\right\}  \right)  $ is the point on $\left\{  i,j\right\}  $ which is a
distance $d_{ij}$ to $\left\{  i\right\}  $ and a distance $d_{ji}$ to
$\left\{  j\right\}  .$ Now, given centers $C\left(  \sigma^{k}\right)  $ for
$k\leq m-1,$ we construct $C\left(  \sigma^{m}\right)  $ as follows. Label the
vertices of $\sigma^{m}$ to be $\left\{  0,1,\ldots,m\right\}  .$

Let $\Pi_{\left\{  i,j\right\}  }$ denote the plane in $\mathbb{R}^{m}$
through $C\left(  \left\{  i,j\right\}  \right)  $ and perpendicular to
$\left\{  i,j\right\}  $ (this is a hyperplane in $\mathbb{R}^{m}$). First we
construct the center of a simplex $\left\{  0,1,2\right\}  $ ($m=2$). One can
embed the simplex in $\mathbb{R}^{2}$ as the three vertices $\left(
0,0\right)  ,$ $\left(  \ell_{01},0\right)  ,$ and $\left(  \ell_{02}%
\cos\gamma_{0},\ell_{02}\sin\gamma_{0}\right)  ,$ where $\gamma_{0}$ is the
angle at vertex $0.$ The centers of the three edges are realized as $C\left(
\left\{  0,1\right\}  \right)  =\left(  d_{01},0\right)  ,$ $C\left(  \left\{
0,2\right\}  \right)  =\left(  d_{02}\cos\gamma_{0},d_{02}\sin\gamma
_{0}\right)  ,$ and $C\left(  \left\{  1,2\right\}  \right)  =\left(
\ell_{01}-d_{12}\cos\gamma_{1},d_{12}\sin\gamma_{1}\right)  .$ Hence
\begin{align*}
\Pi_{\left\{  0,1\right\}  }  &  =\left\{  \left(  d_{01},t\right)
:t\in\mathbb{R}\right\}  ,\\
\Pi_{\left\{  0,2\right\}  }  &  =\left\{  \left(  d_{02}\cos\gamma_{0}%
+t\sin\gamma_{0},d_{02}\sin\gamma_{0}-t\cos\gamma_{0}\right)  :t\in
\mathbb{R}\right\}  ,\\
\Pi_{\left\{  1,2\right\}  }  &  =\left\{  \left(  \ell_{01}-d_{12}\cos
\gamma_{1}+t\sin\gamma_{1},d_{12}\sin\gamma_{1}+t\cos\gamma_{1}\right)
:t\in\mathbb{R}\right\}  .
\end{align*}
A quick calculation (using the law of cosines to compute $\cos\gamma_{i}$ and
$\sin\gamma_{i}$ in terms of $d_{ij}$) shows that the three intersection
points of these lines coincide if and only if
(\ref{compatibility for local duality}) holds.

We now construct $C\left(  \sigma^{m}\right)  $ given $C\left(  \sigma
^{m-1}\right)  $ for all $\left(  m-1\right)  $-dimensional simplices. Since
$\sigma^{m}$ is a nondegenerate Euclidean simplex, the planes $\Pi_{\left\{
0,1\right\}  },\ldots,\Pi_{\left\{  0,m\right\}  }$ intersect at one point,
$c.$ We need only show that the planes $\Pi_{\left\{  i,j\right\}  }$ also
intersect $c.$ This is true because inside $\left\{  0,i,j\right\}  ,$ the
planes $\Pi_{\left\{  0,i\right\}  }$ and $\Pi_{\left\{  0,j\right\}  }$ meet
each other and the plane $\Pi_{\left\{  i,j\right\}  }$ at $C\left(  \left\{
0,i,j\right\}  \right)  .$ Furthermore, since these planes are all
perpendicular to $\left\{  0,i,j\right\}  ,$ the intersection $\Pi_{\left\{
0,i\right\}  }\cap\Pi_{\left\{  i,j\right\}  }$ is equal to the intersection
$\Pi_{\left\{  0,i\right\}  }\cap\Pi_{\left\{  0,j\right\}  }$ and hence
contains $c.$ We call this point $C\left(  \sigma^{m}\right)  =c.$
\end{proof}

\begin{figure}
[ptb]
\begin{center}
\includegraphics[
natheight=4.854000in,
natwidth=6.248800in,
height=2.3864in,
width=3.0663in
]
{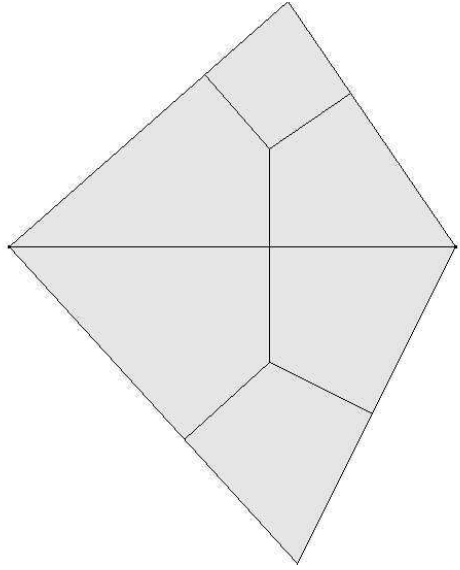}
\caption{Two triangles with the pieces of dual edges intersecting the
triangles included.}
\label{2ddualpic}
\end{center}
\end{figure}

\begin{figure}
[ptbptb]
\begin{center}
\includegraphics[
natheight=7.776400in,
natwidth=8.490400in,
height=2.6263in,
width=2.8652in
]
{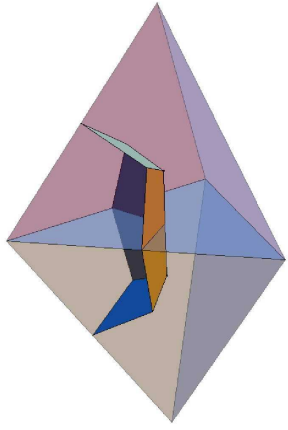}
\caption{Two tetrahedra with the pieces of dual edges and faces intersecting
the tetrahedra included.}
\label{3ddualpic}
\end{center}
\end{figure}

Centers allow a geometric description of the Poincar\'{e} dual of the
triangulation. Any triangulation of a manifold has a cell complex which is its
Poincar\'{e} dual (see, for instance, \cite{Bre} or \cite{Hat}). As noted by
Hirani \cite{Hir}, the assignment of a center to each simplex allows one to
assign a geometric Poincar\'{e} dual, or just dual for short. See Figures
\ref{2ddualpic} and \ref{3ddualpic} for two-dimensional and three-dimensional
simplices with dual cells included. Hirani restricted himself to
\textquotedblleft well-centered\textquotedblright\ triangulations, which means
that the center of each simplex is inside the simplex. This is a very strong
restriction, for even Delaunay triangulations may not be well-centered.
Duality structures allow one to define geometric duals (a realization of the
Poincar\'{e} dual), each of which has a volume. The structure may not be
well-centered, and for this reason some volumes may be negative. The
$k$-dimensional volume of a simplex $\sigma^{k}$ will be denoted $\left\vert
\sigma^{k}\right\vert $ (for instance $\left\vert \left\{  i,j\right\}
\right\vert =\ell_{ij}$) and the $\left(  n-k\right)  $-dimensional (signed)
volume of the dual of a simplex $\bigstar\sigma^{k}$ will be denoted
$\left\vert \bigstar\sigma^{k}\right\vert .$%

It is helpful to consider an example before considering the general
definitions. Given a triangulation of a three-dimensional manifold, one
defines the duals as follows (compare with Figure \ref{3ddualpic}):

\begin{enumerate}
\item[0.] The dual of a 3-simplex $\left\{  i,j,k,\ell\right\}  $ is the
center, $\bigstar\left\{  i,j,k,\ell\right\}  =C\left(  \left\{
i,j,k,\ell\right\}  \right)  ,$ and its volume is one.

\item[1.] The dual of a 2-simplex $\left\{  i,j,k\right\}  $ contained in
$\left\{  i,j,k,\ell\right\}  $ and $\left\{  i,j,k,m\right\}  $ is a 1-cell
$\bigstar\left\{  i,j,k\right\}  ,$ which is the union of the line from
$C\left(  \left\{  i,j,k,\ell\right\}  \right)  $ to $C\left(  \left\{
i,j,k\right\}  \right)  $ and the line from $C\left(  \left\{
i,j,k,m\right\}  \right)  $ to $C\left(  \left\{  i,j,k\right\}  \right)  $.
Its volume is slightly tricky. We define the volume as
\begin{align*}
\left\vert \bigstar\left\{  i,j,k\right\}  \right\vert  &  =\pm d\left[
C\left(  \left\{  i,j,k,\ell\right\}  \right)  ,C\left(  \left\{
i,j,k\right\}  \right)  \right]  \pm d\left[  C\left(  \left\{
i,j,k,m\right\}  \right)  ,C\left(  \left\{  i,j,k\right\}  \right)  \right]
\\
&  =\pm d\left[  C\left(  \left\{  i,j,k,\ell\right\}  \right)  ,C\left(
\left\{  i,j,k,m\right\}  \right)  \right]
\end{align*}
where $d$ is the Euclidean distance in $\mathbb{R}^{3}$ (these are well
defined because we can embed the hinge in $\mathbb{R}^{3}$) and the signs are
defined appropriately. In the first line, the sign is positive if $C\left(
\left\{  i,j,k,\ell\right\}  \right)  $ is on the same side of the plane
containing the side $\left\{  i,j,k\right\}  $ as the simplex $\left\{
i,j,k,\ell\right\}  $ is, and negative if it is on the other side (similarly
for $\left\{  i,j,k,m\right\}  $). The sign on the second line is defined to
be compatible with the previous definition. Note that it is possible for
$\left\vert \bigstar\left\{  i,j,k\right\}  \right\vert $ to be negative.

\item[2.] The dual of a 1-simplex $\left\{  i,j\right\}  $ is the union of
triangles. For each $k,\ell$ such that $\left\{  i,j,k,\ell\right\}  $ is a
simplex, the intersection of the simplex with the dual $\bigstar\left\{
i,j\right\}  $ is the union of the right triangle with vertices $C\left(
\left\{  i,j,k,\ell\right\}  \right)  ,$ $C\left(  \left\{  i,j,k\right\}
\right)  ,$ $C\left(  \left\{  i,j\right\}  \right)  $ and the right triangle
with vertices $C\left(  \left\{  i,j,k,\ell\right\}  \right)  ,$ $C\left(
\left\{  i,j,\ell\right\}  \right)  ,$ $C\left(  \left\{  i,j\right\}
\right)  .$ Each of these triangles has a signed area. The first is
\[
\pm\frac{1}{2}d\left[  C\left(  \left\{  i,j,k,\ell\right\}  \right)
,C\left(  \left\{  i,j,k\right\}  \right)  \right]  ~d\left[  C\left(
\left\{  i,j\right\}  \right)  ,C\left(  \left\{  i,j,k\right\}  \right)
\right]
\]
and the second is defined similarly. The sign is defined as the product of the
appropriate signs in each of the two distances.

\item[3.] The dual of a vertex $\left\{  i\right\}  $ is a union of right
tetrahedra. For each $j,k,\ell$ such that $\left\{  i,j,k,\ell\right\}  $ is a
simplex, the intersection of $\bigstar\left\{  i\right\}  $ with $\left\{
i,j,k,\ell\right\}  $ is the union of six tetrahedra:

\begin{enumerate}
\item the tetrahedron defined by the vertices $C\left(  \left\{
i,j,k,\ell\right\}  \right)  ,$ $C\left(  \left\{  i,j,k\right\}  \right)  ,$
$C\left(  \left\{  i,j\right\}  \right)  ,$ and $i,$

\item the tetrahedron defined by $C\left(  \left\{  i,j,k,\ell\right\}
\right)  ,$ $C\left(  \left\{  i,j,k\right\}  \right)  ,$ $C\left(  \left\{
i,k\right\}  \right)  ,$ and $i,$

\item the tetrahedron defined by $C\left(  \left\{  i,j,k,\ell\right\}
\right)  ,$ $C\left(  \left\{  i,j,\ell\right\}  \right)  ,$ $C\left(
\left\{  i,j\right\}  \right)  ,$ and $i,$

\item the tetrahedron defined by $C\left(  \left\{  i,j,k,\ell\right\}
\right)  ,$ $C\left(  \left\{  i,j,\ell\right\}  \right)  ,$ $C\left(
\left\{  i,\ell\right\}  \right)  ,$ and $i,$

\item the tetrahedron defined by $C\left(  \left\{  i,j,k,\ell\right\}
\right)  ,$ $C\left(  \left\{  i,k,\ell\right\}  \right)  ,$ $C\left(
\left\{  i,k\right\}  \right)  ,$ and $i,$

\item and the tetrahedron defined by $C\left(  \left\{  i,j,k,\ell\right\}
\right)  ,$ $C\left(  \left\{  i,k,\ell\right\}  \right)  ,$ $C\left(
\left\{  i,\ell\right\}  \right)  ,$ and $i.$
\end{enumerate}

\noindent The volume of $\bigstar\left\{  i\right\}  $ is the sum of the
volumes of these tetrahedra, namely%
\[
\pm\frac{1}{6}d\left[  C\left(  \left\{  i,j,k,\ell\right\}  \right)
,C\left(  \left\{  i,j,k\right\}  \right)  \right]  ~d\left[  C\left(
\left\{  i,j\right\}  \right)  ,C\left(  \left\{  i,j,k\right\}  \right)
\right]  ~d\left[  i,C\left(  \left\{  i,j\right\}  \right)  \right]
\]
for the first and similarly for the others, where the signs are defined appropriately.
\end{enumerate}

We can define the geometric duals in a triangulation of an $n$-dimensional
manifold inductively as follows.

\begin{definition}
Define the dual of $\left\{  0,\ldots,n\right\}  $ to be $\bigstar\left\{
0,\ldots,n\right\}  =C\left(  \left\{  0,\ldots,n\right\}  \right)  ,$ and
$\left\vert \bigstar\left\{  0,\ldots,n\right\}  \right\vert =1.$
\end{definition}

\begin{definition}
The signed distance
\[
d_{\pm}\left[  C\left(  \sigma^{n}\right)  ,C\left(  \sigma^{n-1}\right)
\right]
\]
for $\sigma^{n-1}\subset\sigma^{n}$ is equal to the distance between $C\left(
\sigma^{n}\right)  $ and $C\left(  \sigma^{n-1}\right)  $ in any isometric
embedding $\sigma^{n}\subset\mathbb{R}^{n}$ with the sign positive if
$C\left(  \sigma^{n}\right)  $ is on the same side of the hyperplane defined
by $\sigma^{n-1}\subset\mathbb{R}^{n}$ as $\sigma^{n}$ is, and negative if
$C\left(  \sigma^{n}\right)  $ is on the opposite side.
\end{definition}

It will be useful to know the following formula for the distance between the
center of a triangle and the center of a side. Consider a triangle $\left\{
i,j,k\right\}  .$ Then some basic Euclidean geometry yields
\begin{equation}
d_{\pm}\left[  C\left(  \left\{  i,j,k\right\}  \right)  ,C\left(  \left\{
i,j\right\}  \right)  \right]  =\frac{d_{ik}-d_{ij}\cos\gamma_{i}}{\sin
\gamma_{i}} \label{center distance}%
\end{equation}
where $\gamma_{i}$ is the angle at vertex $i.$

\begin{proposition}
\label{volume with signed distance}For any $k\geq1,$ the volume of a simplex
$\sigma^{k}$ is%
\begin{equation}
\left\vert \sigma^{k}\right\vert =\frac{1}{k!}\sum_{\sigma^{0}\subset
\cdots\subset\sigma^{k}}\prod\limits_{j=0}^{k-1}d_{\pm}\left[  C\left(
\sigma^{j}\right)  ,C\left(  \sigma^{j+1}\right)  \right]
\label{volume formula for a simplex}%
\end{equation}
where $\sigma^{k}$ is fixed and the sum is over all strings of simplices
contained in $\sigma^{k}.$
\end{proposition}

\begin{proof}
The proof is by induction on $k.$ If $k=1,$ then $\left\vert \left\{
i,j\right\}  \right\vert =d_{ij}+d_{ji}$. Assume
(\ref{volume formula for a simplex}) is true and consider $\sigma^{k+1}.$ Let
the boundary of $\sigma^{k+1}$ be made up of $\sigma_{0}^{k},\ldots
,\sigma_{k+1}^{k}.$ The volume can be computed as
\[
\left\vert \sigma^{k+1}\right\vert =\frac{1}{k+1}\sum_{i=0}^{k+1}d_{\pm
}\left[  C\left(  \sigma_{i}^{k}\right)  ,C\left(  \sigma^{k+1}\right)
\right]  \left\vert \sigma_{i}^{k}\right\vert
\]
where each term in the sum is the volume of the simplex consisting of the
center $C\left(  \sigma^{k+1}\right)  $ union $\sigma_{i}^{k}$ and the signs
for $d_{\pm}$ tell us whether to add the area or subtract the area. It follows
from the inductive hypothesis that
\[
\left\vert \sigma^{k+1}\right\vert =\frac{1}{\left(  k+1\right)  !}%
\sum_{\sigma^{0}\subset\cdots\subset\sigma^{k+1}}\prod\limits_{j=0}^{k}d_{\pm
}\left[  C\left(  \sigma^{j}\right)  ,C\left(  \sigma^{j+1}\right)  \right]
.
\]

\end{proof}

Note that the above argument works for any choice of center $C\left(
\sigma^{k}\right)  \in\mathbb{R}^{k}$ as long as $C\left(  \sigma^{\ell
}\right)  $ are the orthogonal projections onto the subspaces spanned by
$\sigma^{\ell}$ for each subsimplex. The volume of a dual simplex is defined
as follows.

\begin{definition}
\label{definition of dual volume}The volume of a dual simplex $\bigstar
\sigma^{k}$ is defined to be
\begin{equation}
\left\vert \bigstar\sigma^{k}\right\vert =\frac{1}{\left(  n-k\right)  !}%
\sum_{\sigma^{k}\subset\cdots\subset\sigma^{n}}\prod\limits_{j=k}^{n-1}d_{\pm
}\left[  C\left(  \sigma^{j}\right)  ,C\left(  \sigma^{j+1}\right)  \right]
\label{dual volume fmla}%
\end{equation}
where $\sigma^{k}$ is fixed and the sum is over all strings of simplices
containing $\sigma^{k}.$
\end{definition}

Note that the volume is signed (it may be negative). We note that the total
volume is expressible in terms of volumes of the dual simplices.

\begin{proposition}
Given a duality triangulation $\mathcal{T}$ of dimension $n,$ the total volume
is
\begin{equation}
V=\sum_{\sigma^{n}\in\mathcal{T}_{n}}\left\vert \sigma^{n}\right\vert
=\sum_{i\in\mathcal{T}_{0}}\left\vert \bigstar\left\{  i\right\}  \right\vert
. \label{volumes equal}%
\end{equation}

\end{proposition}

\begin{proof}
We know that
\[
\left\vert \bigstar\left\{  i\right\}  \right\vert =\frac{1}{n!}\sum_{\left\{
i\right\}  \subset\cdots\subset\sigma^{n}}\prod\limits_{j=0}^{n-1}d_{\pm
}\left[  C\left(  \sigma^{j}\right)  ,C\left(  \sigma^{j+1}\right)  \right]
\]
by (\ref{dual volume fmla}) and
\[
\left\vert \sigma^{n}\right\vert =\frac{1}{n!}\sum_{\sigma^{0}\subset
\cdots\subset\sigma^{n}}\prod\limits_{j=0}^{n-1}d_{\pm}\left[  C\left(
\sigma^{j}\right)  ,C\left(  \sigma^{j+1}\right)  \right]
\]
by (\ref{volume formula for a simplex}). Hence it is sufficient to show that
\[
\sum_{i\in\mathcal{T}_{0}}\sum_{\left\{  i\right\}  \subset\cdots\subset
\sigma^{n}}%
\]
is a reordering of
\[
\sum_{\sigma^{n}\in\mathcal{T}_{n}}\sum_{\sigma^{0}\subset\cdots\subset
\sigma^{n}}.
\]
Here is one way to see this. Make a graph whose vertices are all simplices of
all dimensions and whose edges connect two simplices if one simplex is in the
boundary of the other. An easy way to draw the graph in the plane is to put
vertices corresponding to $n$-dimensional simplices in a horizontal line on
top, then $\left(  n-1\right)  $-dimensional simplices in a horizontal line
below those, and so on until at the bottom is a horizontal line containing all
of the vertices corresponding to $0$-dimensional simplices in the
triangulation. Now draw the edges, which can only connect a vertex in a row to
a vertex in the row above or below. Now we shall see that both sums are equal
to the sum over all paths between the top and bottom of this graph. We can
count this in two ways, first start at the bottom with each path starting at a
$0$-dimensional simplex, or first start at the top with each path starting at
an $n$-dimensional simplex. These are the two sums.
\end{proof}

\subsection{Equivalence of metric
triangulations\label{equivalence of triangulations}}

We shall now show that weighted triangulations are equivalent to Thurston
triangulations, and that, up to a universal scaling of the weights, both are
almost equivalent to the set of duality triangulations. This is motivated by
the geometric interpretations of the lengths, weights, angles, etc.

First we show the equivalence of weighted triangulations and Thurston triangulations.

\begin{theorem}
There is a bijection between weighted triangulations and Thurston triangulations.
\end{theorem}

\begin{proof}
The definition of Thurston triangulation gives the map to weighted
triangulations, keeping $w_{i}$ the same and assigning
\[
\ell_{ij}=\sqrt{w_{i}+w_{j}-c_{ij}}.
\]
Since we assumed that $w_{i}+w_{j}-c_{ij}>0$, $\ell_{ij}$ must be positive.
Similarly, we can map the other way as
\[
c_{ij}=w_{i}+w_{j}-\ell_{ij}^{2}.
\]
Note that since $\ell_{ij}>0,$ we must have that $w_{i}+w_{j}-c_{ij}>0$.
\end{proof}

Next we map weighted triangulations to duality triangulations. Notice that
there is a one parameter family of deformations of a given weighted
triangulation of a triangle $\left\{  i,j,k\right\}  $ which fix the center
$C\left(  \left\{  i,j,k\right\}  \right)  $. These deformations are given by%
\begin{equation}
w_{i}\rightarrow w_{i}+t \label{deformation}%
\end{equation}
for varying $t.$ We call these \emph{weight scaling deformations}, or just
\emph{weight scalings}.

\begin{theorem}
\label{equivalence of weighted and duality}Weighted triangulations modulo
weight scalings can be mapped injectively into the set of duality
triangulations. It is a bijection if the set of duality triangulations are
required to satisfy
\begin{equation}
\sum_{k=0}^{r}\left(  d_{i_{k}i_{k-1}}^{2}-d_{i_{k-1}i_{k}}^{2}\right)  =0
\label{loop property}%
\end{equation}
for all loops $j=i_{0},i_{1},\ldots,i_{r}=j$, where $\left\{  i_{k}%
,i_{k+1}\right\}  \in\mathcal{T}_{1}.$
\end{theorem}

\begin{proof}
The key observation is that given spheres at the vertices of a simplex with
given radii $\sqrt{w_{i}},$ one can always construct a sphere which is
orthogonal to each of these spheres. The center of that sphere will be the
center of the simplex, and for that reason is often called the
\emph{orthogonal center }\cite{ES}. By the arguments above, we need only
construct the dual for triangles. One can do this very easily by embedding the
circles in a vector space of signature $1,1,1,-1$ as in \cite[40.2]{Ped}.
Given a center, one can draw the lines perpendicular to the sides of the
triangle through the center, and these determine $d_{ij}.$ A careful
calculation yields%
\begin{equation}
d_{ij}=\frac{\ell_{ij}^{2}+w_{i}-w_{j}}{2\ell_{ij}}. \label{d_ij from (w,l)}%
\end{equation}
This is the map to duality triangulations. Note that the condition
(\ref{compatibility for local duality}) is automatically satisfied.

There appears to be more information in weighted triangulations, however,
because the new circle centered at the orthogonal center has a radius, which
can be calculated to be
\begin{align}
r_{ijk}^{2}  &  =d_{ij}^{2}+\left(  \frac{d_{ik}-d_{ij}\cos\gamma_{ijk}}%
{\sin\gamma_{ijk}}\right)  ^{2}-w_{i}\label{r_ijk}\\
&  =\frac{d_{ij}^{2}+d_{ik}^{2}-2d_{ij}d_{ik}\cos\gamma_{ijk}}{\sin^{2}%
\gamma_{ijk}}-w_{i},\nonumber
\end{align}
where $\gamma_{ijk}$ is the angle at vertex $i$ in triangle $\left\{
i,j,k\right\}  .$ Note that $r_{ijk}^{2}=w_{C\left(  \left\{  i,j,k\right\}
\right)  },$ the weight assigned to the center of $\left\{  i,j,k\right\}  .$
The weight scalings allow, for any single triangle $\left\{  i,j,k\right\}  ,$
one to specify the value of $r_{ijk}^{2}$ while fixing the center $C\left(
\left\{  i,j,k\right\}  \right)  .$ Fixing the center means that each would
map to the same duality triangulation. It is easy to see that the formula
(\ref{d_ij from (w,l)}) is unchanged by scaling deformations like
(\ref{deformation}). If one chooses $r_{ijk}$ then the map is unique. Once
this scale is fixed in one triangle, however, the scale is determined on
adjacent triangles, because weights on shared vertices have been fixed, and
the deformation (\ref{deformation}) must be done for all vertices $i$ in the
triangle. Thus there is one free scaling parameter for the whole triangulation
(if it is connected).
\end{proof}

\begin{figure}
[ptb]
\begin{center}
\includegraphics[
natheight=7.776400in,
natwidth=8.490400in,
height=2.5479in,
width=2.7812in
]
{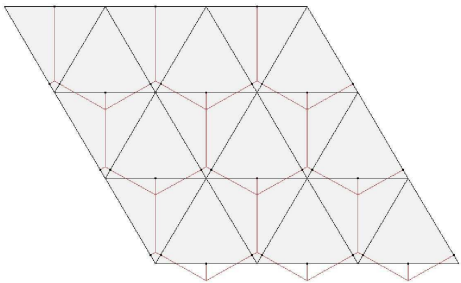}%
\caption{A triangulation of the torus together with dual edges.}
\label{torus}
\end{center}
\end{figure}

The inverse map from duality triangulations to weighted triangulations must
take $d_{ij}+d_{ji}$ to $\ell_{ij}.$ In order to get the weights, we must
first fix $w_{0}$ for a given vertex (this is a free parameter since we are
considering the weighted triangulation modulo scaling). Then each neighboring
weight can be calculated using (\ref{d_ij from (w,l)}):%
\begin{equation}
w_{j}=d_{ji}^{2}-d_{ij}^{2}+w_{i}. \label{w from d_ij}%
\end{equation}
We need only show that this is well defined. Suppose $\left\{  i,j,k\right\}
\in\mathcal{T}_{2}$ and consider a $w_{k}$ which can be defined from $w_{j}$
or $w_{i}.$ Then we need that
\[
d_{ki}^{2}-d_{ik}^{2}+w_{i}=d_{kj}^{2}-d_{jk}^{2}+w_{j}.
\]
But since $w_{j}=d_{ji}^{2}-d_{ij}^{2}+w_{i},$ this follows from the fact that
$d_{ki}^{2}-d_{ik}^{2}=d_{kj}^{2}-d_{jk}^{2}+d_{ji}^{2}-d_{ij}^{2}$ from
(\ref{compatibility for local duality}). It follows by a similar argument that
any null-homotopic loop can be triangulated and property (\ref{loop property})
holds automatically, showing that for any null-homotopic loop $j=i_{0}%
,i_{1},\ldots,i_{L}=j$ of $L$ vertices with $\left\{  i_{k},i_{k+1}\right\}
\in\mathcal{T}_{1}$,
\[
w_{j}=\sum_{k=1}^{L}\left(  d_{i_{k}i_{k-1}}^{2}-d_{i_{k-1}i_{k}}^{2}\right)
+w_{j}.
\]
Thus, in general, we need to assume property (\ref{loop property}) is
satisfied for the weights to be well-defined. For example, the following
triangulation of the torus does not satisfy (\ref{loop property}) for all
loops. Tile a torus with the two triangles $\left\{  1,2,3\right\}  ,\left\{
1,2,4\right\}  $ where $d_{31}=d_{21}=d_{24}=1-\varepsilon,$ $d_{13}%
=d_{12}=d_{42}=\varepsilon,$ and $d_{32}=d_{23}=d_{14}=d_{41}=\frac{1}{2}$ for
small $\varepsilon,$ see Figure \ref{torus}. Note that
\begin{align*}
d_{12}^{2}+d_{23}^{2}+d_{31}^{2}  &  =\varepsilon^{2}+\frac{1}{4}+\left(
1-\varepsilon\right)  ^{2}=d_{21}^{2}+d_{13}^{2}+d_{32}^{2}\\
d_{12}^{2}+d_{24}^{2}+d_{41}^{2}  &  =\varepsilon^{2}+\frac{1}{4}+\left(
1-\varepsilon\right)  ^{2}=d_{21}^{2}+d_{14}^{2}+d_{42}^{2}%
\end{align*}
and so on. The homotopy-nontrivial loop containing $\left\{  1,2\right\}  $
will not satisfy property (\ref{loop property}). However, if we started with a
weighted triangulation, property (\ref{loop property})\ is automatically
satisfied and thus the map from weighted triangulations to duality
triangulations is injective.%

\begin{corollary}
For a triangulation of a simply connected manifold, there is a bijection
between weighted triangulations up to scaling and duality triangulations.
\end{corollary}

\begin{proof}
Since the manifold is simply connected, any loop bounds a 2-dimensional disk,
homeomorphic to $D^{2}=\left\{  x\in\mathbb{R}^{2}:\left\vert x\right\vert
^{2}\leq1\right\}  $, which is triangulated. One can easily prove by induction
on the number of triangles triangulating the disk that on the boundary of any
such disk, (\ref{loop property}) holds.
\end{proof}

\section{Regular triangulations\label{regular triangulations}}

\subsection{Introduction to regular triangulations}

Recall the definition of a regular triangulation (see, for instance, \cite{ES}
or \cite{AK}). Let $d\left(  x,p\right)  $ be the Euclidean distance between
points $p$ and $x.$ Define the power distance
\[
\pi_{p}:\mathbb{R}^{n}\rightarrow\mathbb{R}%
\]
by%
\begin{equation}
\pi_{p}\left(  x\right)  =d\left(  x,p\right)  ^{2}-w_{p}
\label{power definition}%
\end{equation}
if $p$ is a point weighted with $w_{p}.$ The power is important as a function
which is zero on the sphere centered at $p$ with radius $\sqrt{w_{p}}$,
positive outside the sphere, and negative inside the sphere. Notice that if
$p$ is a vertex of a simplex $\sigma$ and $c=C\left(  \sigma\right)  $ then
$\pi_{c}\left(  p\right)  =w_{p}$ and $\pi_{p}\left(  c\right)  =w_{c},$ where
the weight $w_{c}$ is defined as the square of the radius of the orthogonal
sphere as described in Section \ref{weighted triangulations}.

Since we can embed any hinge in $\mathbb{R}^{n},$ the following local
definition of regularity makes sense on a piecewise Euclidean manifold.

\begin{definition}
An $\left(  n-1\right)  $-dimensional simplex $\sigma^{n-1}$ incident on two
$n$-dimensional simplices $\sigma_{1}^{n}=\sigma^{n-1}\cup\left\{
v_{1}\right\}  $ and $\sigma_{2}^{n}=\sigma^{n-1}\cup\left\{  v_{2}\right\}  $
is \emph{locally regular} if $\pi_{c_{1}}\left(  v_{2}\right)  >w_{v_{2}}$ and
$\pi_{c_{2}}\left(  v_{1}\right)  >w_{v_{1}},$ where $c_{i}=C\left(
\sigma_{i}^{n}\right)  $ is the center of $\sigma_{i}^{n}$ for $i=1$ or $2.$
If the weights are all equal to zero, a locally regular simplex is said to be
\emph{locally Delaunay}.
\end{definition}

Sometimes we will instead say that the \emph{hinge} is locally regular. A
hinge is locally Delaunay if and only if it satisfies the local empty
circumsphere property: the sphere circumscribing $\sigma_{1}^{n}$ does not
contain $v_{2}$. This is simply the interpretation of the definition when the
weights are equal to zero. Note that the condition for being locally regular
is unchanged by a weight scaling of the type (\ref{deformation}) due to the
formula (\ref{r_ijk}) for $w_{C\left(  \left\{  i,j,k\right\}  \right)  }.$

There are actually global definitions of regular and Delaunay, since the
definition of power (\ref{power definition}) makes sense globally using the
intrinsic distance (\ref{intrinsic distance}) described in Section
\ref{Euclidean structures basic}.

\begin{definition}
An $n$-dimensional weighted triangulation is \emph{regular} if for every
$\sigma^{n}\in\mathcal{T}_{n},$ we have $\pi_{C\left(  \sigma^{n}\right)
}\left(  v\right)  >w_{v}$ for every vertex $v$ in the complement of
$\sigma^{n}.$ In the case that the weights are all zero, we say the
triangulation is \emph{Delaunay}.
\end{definition}

In the case of two-dimensional Delaunay, the condition on the power says that
for every circle containing at least three vertices, there is no vertex inside
that circle. It is a well known fact that for $n$-dimensional regular
triangulations of points in $\mathbb{R}^{n}$ \cite{AK} and for 2-dimensional
piecewise Euclidean surfaces with zero weights \cite{BS} \cite{Lei} that every
hinge being locally regular is equivalent to the triangulation being regular.
It is likely that the proof in \cite[Chapter 3]{Lei} can be generalized to
regular triangulations of any dimension, but we do not do that here.

The argument in \cite{AK} uses the fact that a geodesic must be a straight
line, and along a geodesic line the power increases in the manner listed
below. To generalize that argument, one needs the following assumption:

\begin{criterion}
\label{monotonicity}Suppose the hinge $\left\{  \sigma_{1}^{n},\sigma_{2}%
^{n},\sigma^{n-1}\right\}  $ is locally regular. Consider a minimizing
geodesic ray $\gamma$ starting at $X_{0}$ which intersects a hinge $\left\{
\sigma_{1}^{n},\sigma_{2}^{n},\sigma^{n-1}\right\}  $ by first entering
$\sigma_{1}^{n}$ and then $\sigma_{2}^{n}.$ The simplex $\sigma^{n-1}$
determines a plane which separates $\sigma_{1}^{n}$ and $\sigma_{2}^{n}$ and
contains all points $x$ such that $\pi_{C\left(  \sigma_{1}^{n}\right)
}\left(  x\right)  =\pi_{C\left(  \sigma_{2}^{n}\right)  }\left(  x\right)  .$
Then $\pi_{C\left(  \sigma_{1}^{n}\right)  }\left(  X_{0}\right)
<\pi_{C\left(  \sigma_{2}^{n}\right)  }\left(  X_{0}\right)  .$
\end{criterion}

One might try to prove Criterion \ref{monotonicity} by \textquotedblleft
developing the geodesic\textquotedblright\ in the plane in the following way
(we consider two dimensions for simplicity). Start with a triangle and embed
it in $\mathbb{R}^{2}.$ For each new triangle which the geodesic goes through,
embed a copy in $\mathbb{R}^{2}$ adjacent to the previous triangle so that it
looks like we are unfolding the manifold. The geodesic must be a straight line
if it does not go through a vertex and so we may try to make comparisons on
this development. Note also that by the following theorem of Gluck, every two
points have a minimizing geodesic between them.

\begin{theorem}
[{\cite[Prop. 2.1]{Sto}}]If a piecewise Euclidean manifold is complete with
respect to the intrinsic distance, in particular if $M$ is a finite
triangulation, then there is at least one minimizing geodesic between any two
points of $M.$
\end{theorem}

The problem with this is that geodesics do go through vertices and even by
varying the endpoints slightly, a minimizing geodesic may still go through the
vertex (see \cite[Figure 14]{MP}). Hence it is not at all clear that Criterion
\ref{monotonicity} is always satisfied.

Note that Bobenko and Springborn \cite{BS} are able to prove that Delaunay is
the same as all edges being locally Delaunay in general by developing the
triangulation (not along a geodesic). Their argument appears to strongly use
the fact that the edges are locally Delaunay (with all weights equal to zero),
but does not use Criterion \ref{monotonicity}.

For completeness, we include the proof for regular triangulations of
$n$-dimensional manifolds, assuming Criterion \ref{monotonicity}, which is
proven using a similar method.

\begin{theorem}
Under the assumptiong of Criterion \ref{monotonicity}, an $n$-dimensional
weighted triangulation is regular if and only if all of its hinges are locally regular.
\end{theorem}

\begin{proof}
This proof is essentially the one seen in \cite{AK} for Delaunay
triangulations. Clearly if the triangulation is regular, then all hinges are
locally regular. Now suppose all of the hinges of a weighted triangulation are
locally regular. Given a vertex $v$ and a simplex $\sigma^{n}$ such that $v$
is not in $\sigma^{n},$ we may consider the line $L$ from $v$ to a point in
the simplex $\sigma^{n}.$ Possibly by adjusting the line slightly, it must
intersect, in order, a sequence of $n$-dimensional simplices $\sigma_{1}%
^{n},\ldots\sigma_{k}^{n}=\sigma^{n}$ where $v$ is in a simplex bordering
$\sigma_{1}^{n}.$ By Criterion \ref{monotonicity} we know that
\[
\pi_{C\left(  \sigma_{i}^{n}\right)  }\left(  v\right)  <\pi_{C\left(
\sigma_{i+1}^{n}\right)  }\left(  v\right)
\]
for $i=1,\ldots,k-1.$ Since the triangulation is locally regular,
\[
w_{v}<\pi_{C\left(  \sigma_{1}^{n}\right)  }\left(  v\right)  .
\]
Stringing these together, we get that
\[
w_{v}<\pi_{C\left(  \sigma^{n}\right)  }\left(  v\right)  .
\]

\end{proof}

Although we have not proven that regular triangulations and locally regular
triangulations are the same, we will often suppress the word \textquotedblleft
local\textquotedblright\ in the rest of this paper, always considering the
local property.

\subsection{Regular triangulations and duality
structures\label{regular triangulations and duality}}

In order to have a definition of locally regular in terms of duality
structures, we first look at the two-dimensional case. A regular hinge
$\left\{  \left\{  i,j,k\right\}  ,\left\{  i,j,\ell\right\}  \right\}  $ must
satisfy
\begin{align*}
\pi_{C\left(  \left\{  i,j,k\right\}  \right)  }\left(  \ell\right)   &
=d\left(  C\left(  \left\{  i,j,k\right\}  \right)  ,\left\{  \ell\right\}
\right)  ^{2}-r_{ijk}^{2}>w_{\ell}\\
\pi_{C\left(  \left\{  i,j,\ell\right\}  \right)  }\left(  k\right)   &
=d\left(  C\left(  \left\{  i,j,\ell\right\}  \right)  ,\left\{  k\right\}
\right)  ^{2}-r_{ij\ell}^{2}>w_{k}.
\end{align*}

\begin{proposition}
The center $C\left(  \left\{  i,j,k\right\}  \right)  $ and radius $r_{ijk}$
are uniquely determined by the three equations%
\begin{align*}
d\left(  C\left(  \left\{  i,j,k\right\}  \right)  ,\left\{  i\right\}
\right)  ^{2}-r_{ijk}^{2}  &  =w_{i}\\
d\left(  C\left(  \left\{  i,j,k\right\}  \right)  ,\left\{  j\right\}
\right)  ^{2}-r_{ijk}^{2}  &  =w_{j}\\
d\left(  C\left(  \left\{  i,j,k\right\}  \right)  ,\left\{  k\right\}
\right)  ^{2}-r_{ijk}^{2}  &  =w_{k}.
\end{align*}

\end{proposition}

\begin{proof}
Put the triangle in Euclidean space with vertices $v_{i}=\vec{0},v_{j},v_{k}.$
We know that $C\left(  \left\{  i,j,k\right\}  \right)  =xv_{j}+yv_{k}$ for
some $x$ and $y$ and let $z$ be the unknown radius. Now we can write the first
two equations as
\begin{align*}
\left\vert xv_{j}+yv_{k}\right\vert ^{2}-z^{2}  &  =w_{i}\\
\left\vert \left(  xv_{j}+yv_{k}\right)  -v_{j}\right\vert ^{2}-z^{2}  &
=w_{j}%
\end{align*}
so%
\[
w_{i}-2v_{j}\cdot\left(  xv_{j}+yv_{k}\right)  +\ell_{ij}^{2}=w_{j}%
\]
which is linear in $x,y.$ Similarly, we have
\[
w_{i}-2v_{k}\cdot\left(  xv_{j}+yv_{k}\right)  +\ell_{ik}^{2}=w_{k}.
\]
So the problem reduces to a linear system
\begin{align*}
w_{i}+\ell_{ij}^{2}-w_{j}  &  =2\ell_{ij}^{2}x+2\ell_{ij}\ell_{ik}\left(
\cos\gamma_{i}\right)  y\\
w_{i}+\ell_{ik}^{2}-w_{k}  &  =2\ell_{ij}\ell_{ik}\left(  \cos\gamma
_{i}\right)  x+2\ell_{ik}^{2}y,
\end{align*}
where $\gamma_{i}$ is the angle at vertex $i,$ with solutions
\begin{align*}
x  &  =\frac{\left(  w_{i}+\ell_{ij}^{2}-w_{j}\right)  \ell_{ik}-\left(
w_{i}+\ell_{ik}^{2}-w_{k}\right)  \ell_{ij}\cos\gamma_{i}}{2\left(  \sin
^{2}\gamma_{i}\right)  \ell_{ij}^{2}\ell_{ik}}\\
y  &  =\frac{\left(  w_{i}+\ell_{ik}^{2}-w_{k}\right)  \ell_{ij}-\left(
w_{i}+\ell_{ij}^{2}-w_{j}\right)  \ell_{ik}\cos\gamma_{i}}{2\left(  \sin
^{2}\gamma_{i}\right)  \ell_{ij}\ell_{ik}^{2}}%
\end{align*}
$\allowbreak$\bigskip and
\[
z^{2}=x^{2}\ell_{ij}^{2}+y^{2}\ell_{ik}^{2}+2xy\ell_{ij}\ell_{ik}\cos
\gamma_{i}-w_{i}.
\]

\end{proof}

\begin{corollary}
\label{boundary regular}If an edge is on the boundary of regular, i.e.
\[
\pi_{C\left(  \left\{  i,j,k\right\}  \right)  }\left(  \ell\right)  =d\left(
C\left(  \left\{  i,j,k\right\}  \right)  ,\left\{  \ell\right\}  \right)
^{2}-r_{ijk}^{2}=w_{\ell},
\]
then $C\left(  \left\{  i,j,k\right\}  \right)  =C\left(  \left\{
i,j,\ell\right\}  \right)  $ and $r_{ijk}=r_{ij\ell}.$
\end{corollary}

\begin{proof}
If $d\left(  C\left(  \left\{  i,j,k\right\}  \right)  ,\ell\right)
^{2}-r_{ijk}^{2}=w_{\ell}$ then $\left(  C\left(  \left\{  i,j,k\right\}
\right)  ,r_{ijk}\right)  $ satisfy the same three equations as $\left(
C\left(  \left\{  i,j,\ell\right\}  \right)  ,r_{ij\ell}\right)  ,$ which
determine these uniquely. Hence they must be equal.
\end{proof}

\begin{corollary}
\label{edge regular duality}An edge $\left\{  i,j\right\}  $ is regular if and
only if $\left\vert \bigstar\left\{  i,j\right\}  \right\vert >0.$
\end{corollary}

\begin{proof}
Clearly $\left\vert \bigstar\left\{  i,j\right\}  \right\vert =0$ on the
boundary of regular as in Corollary \ref{boundary regular} since the centers
are the same. It is clear that $\left\vert \bigstar\left\{  i,j\right\}
\right\vert >0$ if the edge is regular.
\end{proof}

One can now address the case of $n$ dimensions. The corresponding proofs go
through essentially untouched, and one has the following characterization of
regular triangulations.

\begin{proposition}
An $\left(  n-1\right)  $-dimensional simplex $\sigma^{n-1}$ which forms a
hinge with simplices $\sigma_{i}^{n}=\sigma^{n-1}\cup\left\{  i\right\}  $ and
$\sigma_{j}^{n}=\sigma^{n-1}\cup\left\{  j\right\}  $ is regular if and only
if $\left\vert \bigstar\sigma^{n-1}\right\vert >0.$
\end{proposition}

Note that $\bigstar\sigma^{n-1}$ is a one-dimensional simplex, so the property
of being regular has to do with lengths dual to $\left(  n-1\right)
$-simplices being positive. The previous discussion motivates the following
definitions which, in light of Theorem
\ref{equivalence of weighted and duality}, are slight generalizations of those
for weighted triangulations.

\begin{definition}
An $n$-dimensional hinge at simplex $\sigma^{n-1}$ is said to be \emph{locally
regular} if $\left\vert \bigstar\sigma^{n-1}\right\vert >0.$ An $n$%
-dimensional duality triangulation $\mathcal{T}$ is said to be \emph{locally
regular} if $\left\vert \bigstar\sigma^{n-1}\right\vert >0$ for all
$\sigma^{n-1}\in\mathcal{T}_{n-1}.$
\end{definition}

The duality structure is called a Voronoi diagram in the case the
triangulation is Delaunay. Voronoi diagrams can be described in a more direct
way. A point $x$ is in the Voronoi cell $\bigstar\left\{  i\right\}  $ if it
is closer to $i$ than to any other vertex. The boundary of the Voronoi cells
forms the $\left(  n-1\right)  $-dimensional complex called the Voronoi
diagram. The analogue for regular triangulations is called a power diagram. A
point $x$ is in the power cell $\bigstar\left\{  i\right\}  $ if its power
distance $\pi_{i}\left(  x\right)  $ is less than $\pi_{j}\left(  x\right)  $
for any $j\neq i$ (see \cite{AK} \cite{ES}). In the case of regular
triangulations, the duality described in Section \ref{duality triangulations}
is the same as using power diagrams. However, our notion of duality is more
general, making sense for weighted triangulations which are not regular.

An interesting question is how to find a regular triangulation of a given
manifold with given weights. One method of construction is via so called
\textquotedblleft flip algorithms.\textquotedblright\

\subsection{Flips in 2D\label{flips 2d}}

\begin{figure}
[ptb]
\begin{center}
\includegraphics[
natheight=6.616900in,
natwidth=8.490400in,
height=2.3624in,
width=3.0266in
]
{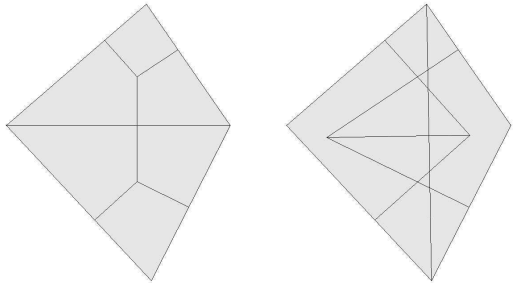}%
\caption{Two hinges differing by a bistellar flip, together with duals.}%
\label{2dflip}%
\end{center}
\end{figure}

We first consider the case of two dimensions. One can imagine the following
notion of a flip. Given a hinge consisting of two triangles $\left\{
i,j,k\right\}  $ and $\left\{  i,j,\ell\right\}  $ incident on one common edge
$\left\{  i,j\right\}  ,$ there exists a flip which exchanges this hinge with
a new hinge, namely $\left\{  i,k,\ell\right\}  $ and $\left\{  j,k,\ell
\right\}  .$ Note that the flip fixes the boundary quadrilateral which
consists cyclically of the vertices $i,k,j,\ell.$ This exchange is called a
$2\rightarrow2$ bistellar flip, or Pachner move (\cite{Pac}). If the hinge is
convex, then this can be done metrically. In fact, the flip can be made at the
level of a duality structure. Given the hinge described above, to do the
bistellar flip we need to construct $d_{k\ell}$ and $d_{\ell k}$ such that the
condition (\ref{compatibility for local duality}) is satisfied in each of the
new triangles. This is done by solving the following system of equations for
$d_{k\ell}$ and $d_{\ell k}$,%
\begin{align*}
d_{ik}^{2}+d_{k\ell}^{2}+d_{\ell i}^{2}  &  =d_{ki}^{2}+d_{i\ell}^{2}+d_{\ell
k}^{2}\\
d_{k\ell}+d_{\ell k}  &  =d\left(  k,\ell\right)
\end{align*}
where $d\left(  k,\ell\right)  $ is the distance between vertex $k$ and vertex
$\ell.$ This distance is the Euclidean distance because the entire hinge can
be embedded in $\mathbb{R}^{2}.$ Note that the first equation is equivalent
to
\[
d_{jk}^{2}+d_{k\ell}^{2}+d_{\ell j}^{2}=d_{kj}^{2}+d_{j\ell}^{2}+d_{\ell
k}^{2}%
\]
using (\ref{compatibility for local duality}) for triangles $\left\{
i,j,k\right\}  $ and $\left\{  i,j,\ell\right\}  .$ The system can actually be
written in a form easier to solve:%
\begin{align}
d_{k\ell}-d_{\ell k}  &  =\frac{d_{ki}^{2}+d_{i\ell}^{2}-d_{\ell i}^{2}%
-d_{ik}^{2}}{d\left(  k,\ell\right)  }\label{flip system}\\
d_{k\ell}+d_{\ell k}  &  =d\left(  k,\ell\right) \nonumber
\end{align}
which is linear, although the dependence of $d\left(  k,\ell\right)  $ on the
remaining $d$'s is not obvious (although easy to find using trigonometry).
Hence the $2\rightarrow2$ bistellar flip is well defined on duality
triangulations, and the triangle inequality follows automatically. The two
hinges which are equivalent by bistellar flips are shown in Figure
\ref{2dflip}.

The flip requires that the quadrilateral is convex, otherwise the flip would
require that one part is folded back, which complicates matters. This
motivates the following definition:

\begin{definition}
A hinge is \emph{flippable }if the quadrilateral defined by the hinge when
embedded in $\mathbb{R}^{2}$ is convex.
\end{definition}

Now, given a convex quadrilateral, there exist two possible ways to make it
into a hinge. The duals are uniquely determined by an assignment of centers to
the edges on the quadrilateral. Let $L_{\left\{  i,j\right\}  }$ be the line
perpendicular to $\left\{  i,j\right\}  $ and through $C\left(  \left\{
i,j\right\}  \right)  .$ Then $L_{\left\{  i,k\right\}  }$ and $L_{\left\{
j,k\right\}  }$ meet at a point which is the center $C\left(  \left\{
i,j,k\right\}  \right)  $ and similarly $L_{\left\{  i,\ell\right\}  }$ and
$L_{\left\{  j,\ell\right\}  }$ meet at a point which is the center $C\left(
\left\{  i,j,\ell\right\}  \right)  .$ However, also $L_{\left\{  i,k\right\}
}$ and $L_{\left\{  i,\ell\right\}  }$ meet at a point which becomes $C\left(
\left\{  i,k,\ell\right\}  \right)  $ after the flip, and similarly with
$L_{\left\{  j,k\right\}  }$ and $L_{\left\{  j,\ell\right\}  }.$ Hence the
centers in the hinge form another quadrilateral dual to the hinge (see the
right side of Figure \ref{2dflip}). One diagonal of the dual quadrilateral
corresponds to $\bigstar\left\{  i,j\right\}  $ and the other corresponds to
$\bigstar\left\{  k,\ell\right\}  .$ One must have positive length and the
other negative length (or both are zero if all dual lines meet at a single
point), so either the hinge is regular, or it will become regular by a flip.
One can also think of the flip of the hinge corresponding to a flip of the
dual hinge. To make this argument rigorous, one simply uses the fact that
$\bigstar\left\{  i,j\right\}  $ must be perpendicular to $\left\{
i,j\right\}  ,$ and considers the possible cases for $\left\vert
\bigstar\left\{  i,j\right\}  \right\vert $ being positive, negative, or zero.
If it is negative, then it must look like the right side Figure \ref{2dflip}
and hence a flip makes $\left\vert \bigstar\left\{  k,\ell\right\}
\right\vert $ positive. If $\left\vert \bigstar\left\{  i,j\right\}
\right\vert $ is zero, then a flip maintains this.

\subsection{Flip algorithms\label{flip algorithms}}

The most naive flip algorithm is to take a given weighted triangulation, look
for a flippable edge which is not regular, and flip it. Continue until the
triangulation is regular. This algorithm was first suggested by Lawson and
shown to find Delaunay triangulations for points in $\mathbb{R}^{2}$
(\cite{Law1}, see also exposition in \cite{Edel} and related result in
\cite{Law2}). It was later shown to work for any 2D piecewise Euclidean
triangulation (where the weights are all zero) independently in \cite{ILTC}
and \cite{Riv}. This turns out not to work to find higher dimensional Delaunay
triangulations or to find regular triangulations (if there are nonzero
weights) even in dimension $2$. It was later found that points in
$\mathbb{R}^{n}$ can be triangulated with regular triangulations (for any
dimension) by incrementally adding one vertex at a time and doing all the
flips before adding additional vertices. In this case one must pay close
attention to the order of the flipping and the algorithm must either sort the
hinges or dynamically decide which hinge to flip next \cite{Joe} \cite{ES}.
Unfortunately, it is not yet clear how to extend these algorithms to piecewise
Euclidean manifolds, since their proofs rely on the fact that the
triangulations are in $\mathbb{R}^{n}.$ In this section we propose a subset of
the space of all weighted triangulations for which the naive flip algorithm
works, just as in the case of two-dimensional Delaunay triangulations.

Consider the following set.

\begin{definition}
A 2-dimensional duality triangulation is said to be \emph{edge positive }if
$d_{ij}>0$ for every directed edge $\left(  i,j\right)  $ of the triangulation
and for any possible flip, i.e. any solution of (\ref{flip system}).
\end{definition}

Hence a triangulation is edge positive if the centers of each edge are inside
the edge and if the center of the new edge after any flip is also inside that
edge. This implies that any non-regular edge is flippable:

\begin{lemma}
\label{2D edge lemma}Given a 2D edge positive duality triangulation, if an
edge is not regular, then it is flippable.
\end{lemma}

\begin{proof}
We prove the contrapositive. Suppose a hinge consisting of $\left\{
i,j,k\right\}  $ and $\left\{  i,j,\ell\right\}  $ is not flippable, i.e. the
quadrilateral is not convex. There can only be one interior angle larger than
$\pi,$ and it must be at vertex $i$ or $j.$ Say it is at $i.$ Let $L_{k}$ be
the line through vertex $i$ which is perpendicular to $\left\{  i,k\right\}  $
and let $L_{\ell}$ be the line through vertex $i$ which is perpendicular to
$\left\{  i,\ell\right\}  $. Since $d_{ik}>0,$ the center $C\left(  \left\{
i,j,k\right\}  \right)  $ must be on the side of $L_{k}$ on which $\left\{
i,k\right\}  $ lies; call this open half-space $H_{k}.$ Similarly, $C\left(
\left\{  i,j,\ell\right\}  \right)  $ must lie on the side of $L_{\ell}$ on
which $\left\{  i,\ell\right\}  $ lies; call this half space $H_{\ell}.$ Let
$H_{j}$ be the half-space containing $\left\{  i,j\right\}  $ whose boundary
is the line $L_{j}$ perpendicular to $\left\{  i,j\right\}  $ through $i.$
Then $C\left(  \left\{  i,j,k\right\}  \right)  $ must be in $H_{k}\cap H_{j}$
and $C\left(  \left\{  i,j,\ell\right\}  \right)  $ must be in $H_{k}\cap
H_{\ell}.$ Since $L_{k},$ $L_{\ell},$ and $L_{j}$ intersect at $i$ and since
the angle at $i$ is more than $\pi,$ $H_{k}\cap H_{j}$ and $H_{\ell}\cap
H_{j}$ are disjoint sectors in a half-space. Use Euclidean isometries to make
put the hinge such that $i$ is at the origin, $\left\{  i,j\right\}  $ is
along the positive $x$-axis, and $k$ has positive $y$-value (and hence $\ell$
must have negative $y$-value). Any possible segment $\bigstar\left\{
i,j\right\}  $ must be on a vertical line which intersects $\left\{
i,j\right\}  .$ It is easy to see that any such line must intersect $H_{k}\cap
H_{j}$ with a larger $y$-value than it intersects $H_{\ell}\cap H_{j},$
implying that $\left\vert \bigstar\left\{  i,j\right\}  \right\vert >0.$
\end{proof}

\begin{theorem}
\label{2D edge flip algo}The edge flip algorithm finds a regular triangulation
given an edge positive duality triangulation.
\end{theorem}

\begin{proof}
Since every flip maintains the edge positive property and every nonregular
edge is flippable, we can always do a flip if the triangulation is not
regular. We now only need an monotone quantity which measures the progress of
the algorithm to complete the proof in the same way as in \cite{AK},
\cite{ES}, \cite{ILTC}, and \cite{Riv}. Since we are in two dimensions, we can
use the Dirichlet energy for almost any function, since the energy increases
if a flip makes the hinge regular (see Theorem \ref{generalized rippa theorem}%
). Since this function increases every time we perform a flip and there are
finitely many possible configurations, the algorithm must terminate.
\end{proof}

Note that the edge flip algorithm to find Delaunay surfaces is a special case,
since in that case, $d_{ij}=\ell_{ij}/2>0.$ In the next section, we suggest
the analogue of this proof for higher dimensions. However, the analogue of
edge positive is possibly less natural in this setting.

\subsection{Higher dimensional flips\label{higher dimensional flips}}

\begin{figure}
[ptb]
\begin{center}
\includegraphics[
natheight=4.702700in,
natwidth=8.506100in,
height=2.51in,
width=4.5275in
]%
{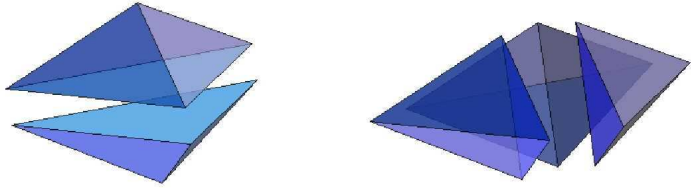}%
\caption{A $2\rightarrow3$ flip. There are two tetrahedra on the left and
three tetrahedra on the right.}%
\label{flip3d}%
\end{center}
\end{figure}

First let's consider the analogue of the $2\rightarrow2$ bistellar move in
higher dimensions. Recall that in any dimension, we can embed a hinge in
$\mathbb{R}^{n},$ so the type of relevant flips must take place inside one or
two simplices in $\mathbb{R}^{n}.$ The relevant flip is the $2\rightarrow n$
flip in $\mathbb{R}^{n}$ (see Figure \ref{flip3d} for the 3D version). The
flip takes two simplices $\sigma_{i}^{n}=\sigma_{0}^{n-1}\cup\left\{
i\right\}  $ and $\sigma_{j}^{n}=\sigma_{0}^{n-1}\cup\left\{  j\right\}  $
meeting at a common face $\sigma_{0}^{n-1}=\left\{  k_{1},\ldots
,k_{n}\right\}  $ and replaces it with $n$ simplices $\sigma_{k_{p}}%
^{n}=\left\{  i,j,k_{1},\ldots,\hat{k}_{p},\ldots,k_{n}\right\}  ,$ where
$\hat{k}_{p}$ indicates that $k_{p}$ is not present. The same argument as
above shows that $d_{ij}$ and $d_{ji}$ can be chosen so that the duality
conditions (\ref{compatibility for local duality}) hold for each face and the
choice is consistent because of the duality conditions which already hold.%

\begin{figure}
[ptb]
\begin{center}
\includegraphics[
natheight=4.702700in,
natwidth=8.506100in,
height=2.689in,
width=4.8503in
]%
{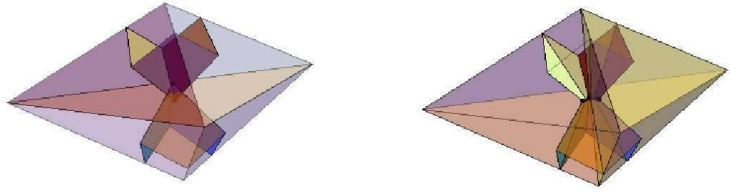}%
\caption{A flip in three dimensions together with dual cells.}%
\label{flip3dwithdual}%
\end{center}
\end{figure}

Now the duality structure gives a hinge a dual hinge similarly to above. Look
at the Figure \ref{flip3dwithdual} to see the 3D case. The boundary of
$\sigma_{i}^{n}$ consists of the faces $\sigma_{0}^{n}=\left\{  k_{1}%
,\ldots,k_{n}\right\}  $ and $\sigma_{ik_{p}}^{n-1}=\left\{  i,k_{1}%
,\ldots,\hat{k}_{p},\ldots,k_{n}\right\}  $ for $p=1,\ldots,n$ while the
boundary of $\sigma_{j}^{n}$ is similarly decomposed. Let $L_{\sigma^{n-1}}$
be the line through $C\left(  \sigma^{n-1}\right)  $ and perpendicular to
$\sigma^{n-1}$ for any $\left(  n-1\right)  $-dimensional simplex. We know
that $L_{\sigma_{ik_{p}}}$ and $L_{\sigma_{ik_{q}}}$ intersect at the point
$C\left(  \sigma_{i}^{n}\right)  $ for every $p,q=1,\ldots,n$ by Proposition
\ref{existence of duals}. We can also consider after the $2\rightarrow n$
flip. The boundary of $\sigma_{k_{p}}^{n}$ consists of $\sigma_{ik_{p}}^{n-1}$
and $\sigma_{jk_{p}}^{n-1}$ together with $\sigma_{k_{p}k_{q}}^{n-1}=\left\{
i,j,k_{1},\ldots,\hat{k}_{p},\ldots,\hat{k}_{q},\ldots,k_{n}\right\}  $ for
$q=1,\ldots,n$ and $q\neq p.$ Hence $L_{\sigma_{ik_{p}}}$ and $L_{\sigma
_{jk_{p}}}$ intersect at the point $C\left(  \sigma_{k_{p}}^{n}\right)  $ for
each $p=1,\ldots,n.$ We find that there is a polytope with vertices $C\left(
\sigma_{i}^{n}\right)  ,$ $C\left(  \sigma_{j}^{n}\right)  ,$ and $C\left(
\sigma_{k_{p}}^{n}\right)  $ for $p=1,\ldots,n.$ This is the dual hinge. The
centers $C\left(  \sigma_{i}^{n}\right)  $ and $C\left(  \sigma_{j}%
^{n}\right)  $ are connected via the edge $\bigstar\sigma_{0}^{n-1}.$ If
$\left\vert \bigstar\sigma_{0}^{n-1}\right\vert <0$ then the flip on the hinge
does a $n\rightarrow2$ flip on the dual hinge which results in removing
$\bigstar\sigma_{0}^{n-1}$ and replaces it with $\bigstar\sigma_{k_{p}k_{q}%
}^{n-1},$ which are $\binom{n}{2}$ dual edges, each with positive length.%

We see that this sort of flipping is exactly what is needed to make regular
triangulations via some sort of flip algorithm. However, the condition of
flippability is harder to guarantee. We now examine flippability.

\begin{definition}
An $n$-dimensional triangulation is said to be $m$-central if $C\left(
\sigma^{k}\right)  $ is inside $\sigma^{k}$ for all $k\leq m.$
\end{definition}

So edge positive is the same as $1$-central. Furthermore, $n$-central is what
is called well-centered in \cite{Hir}. We now show that $\left(  n-1\right)
$-central assures that nonregular hinges are flippable.

\begin{lemma}
Given an $\left(  n-1\right)  $-central triangulation of an $n$-dimensional
manifold, if a hinge is not regular, then it is flippable.
\end{lemma}

\begin{proof}
The proof is essentially the same as the proof of Lemma
\ref{2D edge flip algo}. Consider a hinge consisting of the simplices
$\left\{  i,k_{1},\ldots,k_{n}\right\}  $ and $\left\{  j,k_{1},\ldots
,k_{n}\right\}  $. The first claim is that if the hinge is unflippable, then
at least one dihedral angle must be greater than $\pi.$ This is clear because
if every dihedral angle is less than or equal to $\pi,$ then the hinge is the
intersection of half-spaces defined by the $\left(  n-1\right)  $-simplices on
the boundary and hence convex. Now consider the hyperplanes whose dihedral
angle is greater than $\pi.$ By relabeling we may assume that the hyperplanes
are determined by faces $\sigma_{ik_{n}}^{n-1}=\left\{  i,k_{1},\ldots
,k_{n-1}\right\}  $ and $\sigma_{jk_{n}}^{n-1}=\left\{  j,k_{1},\ldots
,k_{n-1}\right\}  $ and intersect at $\sigma_{0}^{n-2}=\left\{  k_{1}%
,\ldots,k_{n-1}\right\}  .$ Because $C\left(  \sigma_{ik_{n}}^{n-1}\right)
\subset\sigma_{ik_{n}}^{n-1},$ the $C\left(  \sigma_{i}^{n}\right)  $ must be
inside the half-space defined by the plane $\Pi_{ik_{n}}$, the plane through
$\sigma_{0}^{n-2}$ and perpendicular to $\sigma_{ik_{n}}^{n-1},$ on the side
containing $\sigma_{ik_{n}}^{n-1}.$ We have the same for $C\left(  \sigma
_{j}^{n}\right)  $ and since the angle is larger than $\pi$ we must have that
$\left\vert \bigstar\sigma_{0}^{n-1}\right\vert >0$ by a similar argument to
that in the proof of Lemma \ref{2D edge lemma}.
\end{proof}

Regular triangulations of points in $\mathbb{R}^{n}$ are usually produced via
some sort of incremental algorithm (see \cite{ES}, \cite{Joe}). The key
observation is that if a new point is inserted into a regular triangulation,
then there is at least one non-regular hinge which is flippable (or there are
no non-regular hinges and it is regular). The generalization to the manifold
setting is the following. Let $Star\left(  v\right)  ,$ the star of a vertex
$v$, be defined as all simplices containing $v.$

\begin{lemma}
Suppose Criterion \ref{monotonicity} is true. If every hinge in a
triangulation is regular except for hinges intersecting $Star\left(  v\right)
$ for some vertex $v,$ then some if some hinge is not regular, there exists a
flippable nonregular hinge. Hence the triangulation can be made regular via a
flipping algorithm.
\end{lemma}

\begin{proof}
[Proof (sketch)]The proof in \cite{ES} (also with exposition in \cite[Section
12]{Edel}) can be applied to this situation. We are able to prove this lemma
in the generality of manifolds because we have supposed Criterion
\ref{monotonicity} in that generality.
\end{proof}

Using this lemma on subsets of $\mathbb{R}^{n}$, one is able to construct
regular triangulations by: insert one vertex, make the triangulation regular,
and then insert the next vertex, make the triangulation regular, etc.
Unfortunately, on a manifold, it is not clear what the intermediate
triangulations are so the algorithm does not quite work. Also, if one starts
with any triangulation, one may not have a regular triangulation which is
reachable only by flips, as seen in the example \cite[Fig. 5.1]{ES}.

\section{Laplacians\label{laplace}}

Laplace operators on graphs and on piecewise Euclidean manifolds have been
studied in many different contexts, for instance \cite{BS}, \cite{CL},
\cite{Chun}, \cite{G1}, \cite{G2}, \cite{He}, \cite{Hipt}, \cite{Hir},
\cite{MDSB}, \cite{PP}. The purpose of this section is to consider the
comments from Bobenko and Springborn in \cite{BS}, which suggests the use of
Delaunay triangulations as a natural context in which to describe Laplace
operators, and look at the generalization of these comments to regular triangulations.

\subsection{Laplace operator defined\label{laplace defined}}

The suggested Laplace operator on two-dimensional surfaces in \cite{BS} (also
seen in \cite{Hir}, \cite{MDSB}) is the following operator on functions
$f:\mathcal{T}_{0}\rightarrow\mathbb{R},$%
\begin{equation}
\left(  \triangle f\right)  _{i}=\sum_{j:\left\{  i,j\right\}  \in
\mathcal{T}_{1}}w_{ij}\left(  f_{j}-f_{i}\right)  \label{delaunay laplacian}%
\end{equation}
where $w_{ij}$ is defined by
\[
w_{ij}=\frac{1}{2}\left(  \cot\gamma_{kij}+\cot\gamma_{\ell ij}\right)
\]
if $\gamma_{kij}$ is the angle at vertex $k$ in triangle $\left\{
i,j,k\right\}  ,$ and the hinge containing $\left\{  i,j\right\}  $ consists
of the triangles $\left\{  i,j,k\right\}  $ and $\left\{  i,j,\ell\right\}  .$
Note that if $w_{ij}>0$ then this is a Laplacian with weights on the graph
defined by the one-skeleton of the triangulation, and that $\triangle f_{i}>0$
if $f_{i}$ is the minimal value of $f$ and $\triangle f_{i}<0$ if $f_{i}$ is
the maximal value of $f.$ Bobenko and Springborn note that if the
triangulation is Delaunay, then $w_{ij}>0$ and the Laplacian is, in fact, a
Laplacian on graphs in the classical sense (see \cite{Chun}).

A simple calculation shows that if we take the weights at all vertices to be
zero, then the signed distance
\[
d_{\pm}\left[  C\left(  \left\{  i,j,k\right\}  \right)  ,C\left(  \left\{
i,j\right\}  \right)  \right]  =r_{ijk}\cos\gamma_{kij}%
\]
where $r_{ijk}$ is the circumradius of triangle $\left\{  i,j,k\right\}  $.
Since the circumradius can be computed to be
\[
r_{ijk}=\frac{1}{2}\frac{\ell_{ij}}{\sin\gamma_{kij}}%
\]
we find that
\[
d_{\pm}\left[  C\left(  \left\{  i,j,k\right\}  \right)  ,C\left(  \left\{
i,j\right\}  \right)  \right]  =\frac{1}{2}\ell_{ij}\cot\gamma_{kij}.
\]
It immediately follows that%
\[
w_{ij}=\frac{\left\vert \bigstar\left\{  i,j\right\}  \right\vert }{\left\vert
\left\{  i,j\right\}  \right\vert }.
\]
We see that the Delaunay condition is equivalent to $w_{ij}>0,$ which is
equivalent to $\left\vert \bigstar\left\{  i,j\right\}  \right\vert >0.$

In general, Hirani \cite{Hir} suggests the following definition of Laplacian:%
\begin{equation}
\left(  \triangle f\right)  _{i}=\frac{1}{\left\vert \bigstar\left\{
i\right\}  \right\vert }\sum_{j:\left\{  i,j\right\}  \in\mathcal{T}_{1}}%
\frac{\left\vert \bigstar\left\{  i,j\right\}  \right\vert }{\left\vert
\left\{  i,j\right\}  \right\vert }\left(  f_{j}-f_{i}\right)  .
\label{laplacian}%
\end{equation}
This formula has roots in the following integration by parts formula for the
smooth Laplacian:%
\begin{equation}
\int_{U}\triangle f~dV=\int_{\partial U}\nabla f\cdot n~dS
\label{integration by parts}%
\end{equation}
where $n$ is the unit normal to $\partial U.$ Taking $U=\bigstar\left\{
i\right\}  $ and slightly rearranging terms, we get the corresponding formula
on piecewise Euclidean manifolds%
\[
\left(  \triangle f\right)  _{i}~\left\vert \bigstar\left\{  i\right\}
\right\vert =\sum_{j:\left\{  i,j\right\}  \in\mathcal{T}_{1}}\frac
{f_{j}-f_{i}}{\left\vert \left\{  i,j\right\}  \right\vert }\left\vert
\bigstar\left\{  i,j\right\}  \right\vert
\]
where $\frac{f_{j}-f_{i}}{\left\vert \left\{  i,j\right\}  \right\vert }$ is
the normal derivative and $\left\vert \bigstar\left\{  i,j\right\}
\right\vert $ is the surface area measure on the boundary of $\bigstar\left\{
i\right\}  .$ This formula is well defined on any duality triangulation (which
is the motivation for the definition) and coincides with
(\ref{delaunay laplacian}) in the case of Delaunay triangulations, except for
the factor of $\left\vert \bigstar\left\{  i\right\}  \right\vert $. One can
think of the difference between considering the induced measure $\triangle
f~dV$ instead of the pointwise Laplacian $\triangle f.$ It is, in fact,
natural to consider the measure instead since, if we consider the discrete
Laplacian approximating a smooth one, the pointwise Laplacian is only accurate
when considered on scales larger than the scale of the discretization.

We note that the Laplacian given by (\ref{laplacian}) is also the same as the
Laplacian considered by Chow-Luo \cite{CL} in two dimensions as observed by Z.
He, where the duality is defined by Thurston triangulations as described
above. It also appears in \cite{G1} \cite{G2} in three dimensions, where
Thurston triangulations are considered such that $d_{ij}$ depend only on $i.$
Also, the Laplacian described in \cite{Luo} is actually the Laplacian
described above in (\ref{delaunay laplacian}) with the same weights $w_{ij}.$
The interest in these Laplacians is that they are not derived from means such
as (\ref{integration by parts}) but instead as the induced time derivative of
curvature quantities under geometric evolutions.

The Laplacian defined in (\ref{laplacian}) is a Laplacian with weights on
graphs in the usual sense (see \cite{Chun}) if the coefficients
\[
\frac{\left\vert \bigstar\left\{  i,j\right\}  \right\vert }{\left\vert
\bigstar\left\{  i\right\}  \right\vert }%
\]
are each nonnegative. In two dimensions we see that this is implied by
$d_{ij}>0$ and $\left\vert \bigstar\left\{  i,j\right\}  \right\vert >0,$
which is the condition that the triangulation is regular.

Note that the Laplacian can be considered the gradient of a Dirichlet energy
functional as described in \cite{BS}, which is the analogue of the smooth
functional
\[
E\left(  f\right)  =\int_{M}\left\vert \nabla f\right\vert ^{2}dV.
\]
The Dirichlet energy functional induced by the duality triangulation is
\begin{equation}
E\left(  f\right)  =\frac{1}{2}\sum_{\left\{  i,j\right\}  \in\mathcal{T}_{1}%
}\frac{\left\vert \bigstar\left\{  i,j\right\}  \right\vert }{\left\vert
\left\{  i,j\right\}  \right\vert }\left(  f_{j}-f_{i}\right)  ^{2}.
\label{dirichlet energy}%
\end{equation}
This specializes in the case where the $w_{i}=0$ for all $i\in\mathcal{T}_{0}$
(or, equivalently, $d_{ij}=d_{ji}=\ell_{ij}/2$ for all $\left\{  i,j\right\}
\in\mathcal{T}_{1}$) to the Dirichlet energy in \cite{BS}. Note that this
energy is positive if $\left\vert \bigstar\left\{  i,j\right\}  \right\vert
>0.$

\subsection{A generalization of Rippa's theorem\label{rippa}}

Rippa \cite{Rip} showed that if one considers the Dirichlet energy
(\ref{dirichlet energy}) on a triangulation of points in $\mathbb{R}^{2}$
where the weights are zero (or equivalently, $d_{ij}=d_{ji}=\ell_{ij}/2$ for
all edges $\left\{  i,j\right\}  $), flipping to make an edge Delaunay
increases the Dirichlet energy. Bobenko and Springborn \cite{BS} note that his
proof extends trivially to piecewise Euclidean surfaces (2-dimensional
manifolds). We shall express Rippa's theorem in a way closer to the exposition
on \cite{BS}, which is in line with the notation in this paper.

\begin{theorem}
[\cite{Rip}]Let $\left(  \mathcal{T},\ell\right)  $ be a piecewise Euclidean,
triangulated surface with assigned edge lengths $\ell,$ which we think of as a
weighted triangulation with all weights equal to zero. Let $\mathcal{T}_{0}$
be the vertices of the triangulation and let $f:\mathcal{T}_{0}\rightarrow
\mathbb{R}$ be a function. Suppose $\mathcal{T}^{\prime}$ is another
triangulation which is gotten from $\mathcal{T}$ by a $2\rightarrow2$
bistellar flip on edge $e$ (in particular, $\mathcal{T}_{0}=\mathcal{T}%
_{0}^{\prime},$) such that the hinge is locally Delaunay after the flip. Then
\[
E_{\mathcal{T}^{\prime}}\left(  f\right)  \leq E_{\mathcal{T}}\left(
f\right)  ,
\]
where $E_{\mathcal{T}}$ and $E_{\mathcal{T}^{\prime}}$ are the Dirichlet
energies corresponding to $\mathcal{T}$ and $\mathcal{T}^{\prime}.$ As a
consequence, the minimum is attained when all edges are Delaunay (and hence
the triangulation is a Delaunay triangulation).
\end{theorem}

Rippa's proof involves calculating $E\left(  f_{\mathcal{T}^{\prime}}\right)
-E\left(  f_{\mathcal{T}}\right)  $ and showing that it is negative. The key
is a lemma which factors $E\left(  f_{\mathcal{T}^{\prime}}\right)  -E\left(
f_{\mathcal{T}}\right)  $ and for which we shall give a direct proof later for
the more general case of regular triangulations. The only thing missing is the
proof of the final sentence, which requires that flipping edges eventually
produces a Delaunay triangulation, which is proved in \cite{ILTC} and
\cite{Riv}. We can generalize the first part of Rippa's theorem to regular triangulations:

\begin{theorem}
\label{generalized rippa theorem}Let $\left(  \mathcal{T},d\right)  $ be a
duality triangulation of a surface with assigned local lengths $d$. Let
$\mathcal{T}_{0}$ be the vertices of the triangulation and let $f:\mathcal{T}%
_{0}\rightarrow\mathbb{R}$ be a function. Suppose $\left(  \mathcal{T}%
^{\prime},d^{\prime}\right)  $ is another duality triangulation which is
gotten from $\left(  \mathcal{T},d\right)  $ by a $2\rightarrow2$ bistellar
flip on edge $e$ such that the hinge is locally regular after the flip. Then
\[
E_{\mathcal{T}^{\prime}}\left(  f\right)  \leq E_{\mathcal{T}}\left(
f\right)  ,
\]
where $E_{\mathcal{T}}$ and $E_{\mathcal{T}^{\prime}}$ are the Dirichlet
energies corresponding to $\left(  \mathcal{T},d\right)  $ and $\left(
\mathcal{T}^{\prime},d^{\prime}\right)  .$
\end{theorem}

The proof depends on the following important generalization of Rippa's key
lemma \cite[Lemma 2.2]{Rip} (see also \cite{Pow}).

\begin{lemma}
Let $\mathcal{T=}\left\{  \left\{  1,2,3\right\}  ,\left\{  1,2,4\right\}
\right\}  $ and $\mathcal{T}^{\prime}=\left\{  \left\{  1,3,4\right\}
,\left\{  2,3,4\right\}  \right\}  $ be two hinges differing by a flip along
$\left\{  1,2\right\}  $. Then
\[
E\left(  f_{\mathcal{T}^{\prime}}\right)  -E\left(  f_{\mathcal{T}}\right)
=\left(  f_{\mathcal{T}^{\prime}}\left(  c\right)  -f_{\mathcal{T}}\left(
c\right)  \right)  ^{2}A_{1234}^{2}\Phi
\]
where
\[
\Phi=\frac{2\left(  r_{3}r_{4}-r_{1}r_{2}\right)  A_{1234}+w_{1}A_{234}%
+w_{2}A_{134}-w_{3}A_{124}-w_{4}A_{123}}{8A_{123}A_{134}A_{234}A_{124}},
\]
$A_{ijk}$ is the area of $\left\{  i,j,k\right\}  ,$ $A_{1234}=A_{123}%
+A_{124}=A_{134}+A_{234}$ is the area of the hinge, $c$ is the intersection of
the diagonals, $r_{i}$ is the distance between $c$ and vertex $i,$ and
$f_{\mathcal{T}^{\prime}}$ and $f_{\mathcal{T}}$ are the piecewise linear
interpolations of $f$ with respect to the different triangulations. One can
write%
\begin{align*}
f_{\mathcal{T}}\left(  c\right)   &  =\frac{r_{1}}{\ell_{12}}f_{2}+\frac
{r_{2}}{\ell_{12}}f_{1}\\
f_{\mathcal{T}^{\prime}}\left(  c\right)   &  =\frac{r_{3}}{\ell_{34}}%
f_{4}+\frac{r_{4}}{\ell_{34}}f_{3}.
\end{align*}

\end{lemma}

The proof is somewhat involved although straightforward. We use a proof which
is more direct than the ones given by Rippa \cite{Rip} and Powar \cite{Pow}
for the case of Delaunay triangulations.

\begin{proof}
Because we are on a single hinge, it is equivalent to use weighted
triangulations by Theorem \ref{equivalence of weighted and duality}. Let
$\left(  \ell,w\right)  $ be the corresponding lengths and weights. A simple
calculation tells us that
\[
\frac{d_{\pm}\left(  C\left(  \left\{  i,j\right\}  \right)  ,C\left(
\left\{  i,j,k\right\}  \right)  \right)  }{\ell_{ij}}=\frac{1}{2}\cot
\gamma_{kij}+\frac{w_{i}}{2\ell_{ij}^{2}}\cot\gamma_{jik}+\frac{w_{j}}%
{2\ell_{ij}^{2}}\cot\gamma_{ijk}-\frac{w_{k}}{4A_{ijk}},
\]
where $\gamma_{ijk}$ is the angle at vertex $i$ in triangle $\left\{
i,j,k\right\}  $ and $A_{ijk}=\left\vert \left\{  i,j,k\right\}  \right\vert $
is the area. For simplicity, we shall use the notation $h_{ij,k}=d_{\pm
}\left(  C\left(  \left\{  i,j\right\}  \right)  ,C\left(  \left\{
i,j,k\right\}  \right)  \right)  ,$ which we think of as the height of the
triangle $\left\{  i,j,C\left(  \left\{  i,j,k\right\}  \right)  \right\}  .$
Note that $\left\vert \bigstar\left\{  1,2\right\}  \right\vert =h_{12,3}%
+h_{12,4},$ for instance. For any function $f,$ we can compute%
\[
E\left(  f_{\mathcal{T}^{\prime}}\right)  -E\left(  f_{\mathcal{T}}\right)
=\frac{1}{2}\sum_{i,j=1}^{4}a_{ij}f_{i}f_{j},
\]
where
\begin{align*}
a_{12}  &  =\frac{h_{12,3}}{\ell_{12}}+\frac{h_{12,4}}{\ell_{12}}%
,\;\;\;a_{13}=\frac{h_{13,2}}{\ell_{13}}-\frac{h_{13,4}}{\ell_{13}},\\
a_{14}  &  =\frac{h_{14,2}}{\ell_{14}}-\frac{h_{14,3}}{\ell_{14}}%
,\;\;\;a_{23}=\frac{h_{23,1}}{\ell_{23}}-\frac{h_{23,4}}{\ell_{23}},\\
a_{24}  &  =\frac{h_{24,1}}{\ell_{24}}-\frac{h_{24,3}}{\ell_{24}}%
,\;\;\;a_{34}=-\frac{h_{34,1}}{\ell_{34}}-\frac{h_{34,2}}{\ell_{34}},
\end{align*}
and $a_{ii}=-\sum_{j\neq i}a_{ij}$ (where we have symmetrized $a_{ij}=a_{ji}%
$). We now wish to factor the coefficients.

We can easily figure out $r_{i}$ in terms of areas in the following way. For a
realization of the hinge, with $v_{i}$ representing the coordinates of
$\left\{  i\right\}  ,$ we see that $c=v_{1}+\frac{r_{1}}{\ell_{12}}\left(
v_{2}-v_{1}\right)  =v_{3}+\frac{r_{3}}{\ell_{13}}\left(  v_{4}-v_{3}\right)
.$ By taking the cross product with $v_{2}-v_{1}$ or $v_{4}-v_{3}$ we find
that
\[
r_{1}=\frac{\ell_{12}A_{134}}{A_{1234}}~~\text{and~~}r_{3}=\frac{\ell
_{34}A_{123}}{A_{1234}},
\]
where $A_{1234}=A_{123}+A_{124}=A_{134}+A_{234}$ is the area of the entire
hinge. Similarly,
\[
r_{2}=\frac{\ell_{12}A_{234}}{A_{1234}}~~\text{and~~}r_{4}=\frac{\ell
_{34}A_{124}}{A_{1234}}.
\]
Thus
\begin{align*}
f_{\mathcal{T}^{\prime}}\left(  c\right)  -f_{\mathcal{T}}\left(  c\right)
&  =\frac{r_{3}}{\ell_{34}}f_{4}+\frac{r_{4}}{\ell_{34}}f_{3}-\frac{r_{1}%
}{\ell_{12}}f_{2}-\frac{r_{2}}{\ell_{12}}f_{1}\\
&  =\frac{1}{A_{1234}}\left(  A_{123}f_{4}+A_{124}f_{3}-A_{134}f_{2}%
-A_{234}f_{1}\right)  .
\end{align*}
Also useful will be the calculation%
\[
r_{3}r_{4}-r_{1}r_{2}=\frac{1}{A_{1234}^{2}}\left(  \ell_{34}^{2}%
A_{123}A_{124}-\ell_{12}^{2}A_{234}A_{134}\right)  .
\]

There are essentially two different types of coefficients to consider. We need
only consider $a_{12}$ and $a_{13}$ since the others are similar. Let
$\gamma_{ijk}$ be the angle at vertex $i$ in triangle $\left\{  i,j,k\right\}
.$ Consider $a_{12}.$
\begin{align*}
a_{12}  &  =\frac{h_{12,3}}{\ell_{12}}+\frac{h_{12,4}}{\ell_{12}}\\
&  =\frac{1}{2}\cot\gamma_{312}+\frac{w_{1}}{2\ell_{12}^{2}}\cot\gamma
_{213}+\frac{w_{2}}{2\ell_{12}^{2}}\cot\gamma_{123}-\frac{w_{3}}{4A_{123}}\\
&  \;\;\;+\frac{1}{2}\cot\gamma_{412}+\frac{w_{1}}{2\ell_{12}^{2}}\cot
\gamma_{214}+\frac{w_{2}}{2\ell_{12}^{2}}\cot\gamma_{124}-\frac{w_{4}%
}{4A_{124}}\\
&  =\frac{1}{2}\left(  \cot\gamma_{312}+\cot\gamma_{412}\right)  +\frac{w_{1}%
}{2\ell_{12}^{2}}\left(  \cot\gamma_{213}+\cot\gamma_{214}\right) \\
&  \;\;\;+\frac{w_{2}}{2\ell_{12}^{2}}\left(  \cot\gamma_{123}+\cot
\gamma_{124}\right)  -\frac{w_{3}}{4A_{123}}-\frac{w_{4}}{4A_{124}}.
\end{align*}
Let $\theta$ be the angle at $c$ in the triangle $\left\{  1,3,c\right\}  .$
We shall use the fact that in any triangle $\left\{  i,j,k\right\}  $ we have
$\ell_{ij}=\ell_{ik}\cos\gamma_{ijk}+\ell_{jk}\cos\gamma_{jik}$ to compute the
parts.%
\begin{align*}
\cot\gamma_{312}+\cot\gamma_{412}  &  =\frac{\ell_{13}\ell_{23}\cos
\gamma_{312}}{2A_{123}}+\frac{\ell_{14}\ell_{24}\cos\gamma_{412}}{2A_{124}}\\
&  =\frac{\ell_{13}^{2}-\ell_{12}\ell_{13}\cos\gamma_{123}}{2A_{123}}%
+\frac{\ell_{14}^{2}-\ell_{12}\ell_{14}\cos\gamma_{124}}{2A_{124}}\\
&  =\frac{\ell_{13}^{2}}{2A_{123}}+\frac{\ell_{14}^{2}}{2A_{124}}-\left(
\left(  \frac{\sin\gamma_{314}}{\sin\theta\sin\gamma_{123}}-\cot\theta\right)
+\left(  \frac{\sin\gamma_{413}}{\sin\theta\sin\gamma_{124}}+\cot
\theta\right)  \right) \\
&  =\frac{\ell_{13}^{2}}{2A_{123}}+\frac{\ell_{14}^{2}}{2A_{124}}-\frac
{1}{\sin\theta}\left(  \frac{\sin\gamma_{314}}{\sin\gamma_{123}}+\frac
{\sin\gamma_{413}}{\sin\gamma_{124}}\right) \\
&  =\frac{\ell_{13}^{2}}{2A_{123}}+\frac{\ell_{14}^{2}}{2A_{124}}-\frac
{1}{\sin\theta}\frac{\ell_{12}A_{134}A_{1234}}{\ell_{34}A_{123}A_{124}}\\
&  =\frac{\ell_{13}^{2}A_{124}+\ell_{14}^{2}A_{123}-\ell_{12}^{2}A_{134}%
}{2A_{123}A_{124}}\\
&  =\frac{\ell_{13}^{2}+\ell_{14}^{2}}{2A_{1234}}+\frac{\ell_{13}^{2}%
A_{124}^{2}+\ell_{14}^{2}A_{123}^{2}-\ell_{12}^{2}A_{134}^{2}}{2A_{123}%
A_{124}A_{1234}}-\frac{\ell_{12}^{2}A_{134}A_{234}}{2A_{123}A_{124}A_{1234}}%
\end{align*}
since
\[
\sin\gamma_{314}=\cos\gamma_{123}\sin\theta+\sin\gamma_{123}\cos\theta
\]
and
\[
\sin\gamma_{413}=\cos\gamma_{124}\sin\theta-\sin\gamma_{124}\cos\theta.
\]
Furthermore,%
\begin{align*}
\ell_{13}^{2}A_{124}^{2}+\ell_{14}^{2}A_{123}^{2}-\ell_{12}^{2}A_{134}^{2}  &
=\frac{1}{4}\ell_{12}^{2}\ell_{13}^{2}\ell_{14}^{2}\left(  \sin^{2}%
\gamma_{124}+\sin^{2}\gamma_{123}-\sin^{2}\left(  \gamma_{123}+\gamma
_{124}\right)  \right) \\
&  =-\frac{1}{2}\ell_{12}^{2}\ell_{13}^{2}\ell_{14}^{2}\left(  \sin
\gamma_{123}\sin\gamma_{124}\cos\gamma_{134}\right) \\
&  =-2A_{123}A_{124}\ell_{13}\ell_{14}\cos\gamma_{134}%
\end{align*}
since
\[
\sin^{2}A+\sin^{2}B-\sin^{2}\left(  A+B\right)  =-2\sin A\sin B\cos\left(
A+B\right)  .
\]
Thus we have%
\begin{align*}
\cot\gamma_{312}+\cot\gamma_{412}  &  =\frac{\left(  \ell_{13}^{2}+\ell
_{14}^{2}-2\ell_{13}\ell_{14}\cos\gamma_{134}\right)  }{2A_{1234}}-\frac
{\ell_{12}^{2}A_{134}A_{234}}{2A_{123}A_{124}A_{1234}}\\
&  =\frac{\ell_{34}^{2}A_{123}A_{124}-\ell_{12}^{2}A_{134}A_{234}}%
{2A_{1234}A_{123}A_{124}}\\
&  =\frac{A_{1234}}{A_{123}A_{124}}\left(  r_{3}r_{4}-r_{1}r_{2}\right)  .
\end{align*}
For the other parts,%
\begin{align*}
\cot\gamma_{213}+\cot\gamma_{214}  &  =\frac{\cos\gamma_{213}}{\sin
\gamma_{213}}+\frac{\cos\gamma_{214}}{\sin\gamma_{214}}\\
&  =\frac{\sin\gamma_{234}}{\sin\gamma_{213}\sin\gamma_{214}}\\
&  =\frac{\ell_{12}^{2}A_{234}}{2A_{123}A_{124}}%
\end{align*}
and%
\[
\cot\gamma_{123}+\cot\gamma_{124}=\frac{\ell_{12}^{2}A_{134}}{2A_{123}A_{124}%
}.
\]
Thus%
\begin{align*}
a_{12}  &  =\frac{1}{2}\left(  \cot\gamma_{312}+\cot\gamma_{412}\right)
+\frac{w_{1}}{2\ell_{12}^{2}}\left(  \cot\gamma_{213}+\cot\gamma_{214}\right)
\\
&  +\frac{w_{2}}{2\ell_{12}^{2}}\left(  \cot\gamma_{123}+\cot\gamma
_{124}\right)  -\frac{w_{3}}{4A_{123}}-\frac{w_{4}}{4A_{124}}.
\end{align*}
implies that%
\begin{align*}
a_{12}  &  =\frac{A_{234}A_{134}}{4A_{123}A_{134}A_{234}A_{124}}\left(
2A_{1234}\left(  r_{3}r_{4}-r_{1}r_{2}\right)  +w_{1}A_{234}+w_{2}%
A_{134}-w_{3}A_{124}-w_{4}A_{123}\right) \\
&  =2A_{234}A_{134}\Phi.
\end{align*}
(Recall
\[
\Phi=\frac{2A_{1234}\left(  r_{3}r_{4}-r_{1}r_{2}\right)  +w_{1}A_{234}%
+w_{2}A_{134}-w_{3}A_{124}-w_{4}A_{123}}{8A_{123}A_{134}A_{234}A_{124}}%
\]
as in the statement of the lemma.)

Now consider $a_{13}.$ We can compute
\begin{align*}
a_{13}  &  =\frac{h_{13,2}}{\ell_{13}}-\frac{h_{13,4}}{\ell_{13}}\\
&  =\frac{1}{2}\cot\gamma_{213}+\frac{w_{1}}{2\ell_{13}^{2}}\cot\gamma
_{312}+\frac{w_{3}}{2\ell_{13}^{2}}\cot\gamma_{123}-\frac{w_{2}}{4A_{123}}\\
&  \;\;\;-\left(  \frac{1}{2}\cot\gamma_{413}+\frac{w_{1}}{2\ell_{13}^{2}}%
\cot\gamma_{314}+\frac{w_{3}}{2\ell_{13}^{2}}\cot\gamma_{134}-\frac{w_{4}%
}{4A_{134}}\right) \\
&  =\frac{1}{2}\left(  \cot\gamma_{213}-\cot\gamma_{413}\right)  +\frac{w_{1}%
}{2\ell_{13}^{2}}\left(  \cot\gamma_{312}-\cot\gamma_{314}\right) \\
&  \;\;\;+\frac{w_{3}}{2\ell_{13}^{2}}\left(  \cot\gamma_{123}-\cot
\gamma_{134}\right)  -\frac{w_{2}}{4A_{123}}+\frac{w_{4}}{4A_{134}}.
\end{align*}
We see that
\begin{align*}
\cot\gamma_{213}-\cot\gamma_{413}  &  =\frac{\sin\gamma_{324}}{\sin
\gamma_{213}\sin\theta}-\frac{\sin\gamma_{124}}{\sin\gamma_{413}\sin\theta}\\
&  =\frac{\ell_{12}^{2}A_{134}A_{234}-\ell_{34}^{2}A_{123}A_{124}}%
{2A_{1234}A_{123}A_{134}}%
\end{align*}
since $\sin\gamma_{324}=-\cos\theta\sin\gamma_{213}+\sin\theta\cos\gamma
_{213}$ and similarly $\sin\gamma_{124}=-\cos\theta\sin\gamma_{413}+\sin
\theta\cos\gamma_{413}.$ We also get%
\begin{align*}
\cot\gamma_{312}-\cot\gamma_{314}  &  =\frac{\cos\gamma_{312}\sin\gamma
_{314}-\cos\gamma_{314}\sin\gamma_{312}}{\sin\gamma_{312}\sin\gamma_{314}}\\
&  =-\frac{\sin\gamma_{324}}{\sin\gamma_{312}\sin\gamma_{314}}\\
&  =-\frac{\ell_{13}^{2}A_{234}}{2A_{123}A_{134}}%
\end{align*}
and%
\[
\cot\gamma_{123}-\cot\gamma_{134}=\frac{\ell_{13}^{2}A_{124}}{2A_{123}A_{134}%
}.
\]
And so
\begin{align*}
a_{13}  &  =\frac{-A_{234}A_{124}}{4A_{123}A_{134}A_{234}A_{124}}\left(
2A_{1234}\left(  r_{3}r_{4}-r_{1}r_{2}\right)  +w_{1}A_{234}-w_{3}%
A_{124}+w_{2}A_{134}-w_{4}A_{123}\right) \\
&  =-2A_{234}A_{124}\Phi.
\end{align*}

A similar argument gives the other coefficients. Then we see, for instance,
that
\begin{align*}
a_{11}  &  =-a_{12}-a_{13}-a_{14}\\
&  =2\left(  -A_{234}A_{134}+A_{234}A_{124}+A_{234}A_{123}\right)  \Phi\\
&  =2A_{234}^{2}\Phi
\end{align*}
with similar expressions for $a_{22},$ $a_{33},$ and $a_{44}.$ Finally, we get
that
\[
E\left(  f_{\mathcal{T}^{\prime}}\right)  -E\left(  f_{\mathcal{T}}\right)
=\left(  A_{123}f_{4}+A_{124}f_{3}-A_{234}f_{1}-A_{134}f_{2}\right)  ^{2}%
\Phi,
\]
which is equivalent to the lemma.
\end{proof}

Now we can prove the theorem.

\begin{proof}
[Proof of Theorem \ref{generalized rippa theorem}]Note that since we are only
concerned with a hinge, it is equivalent to consider weighted triangulations
or duality triangulations. Since the coefficient $a_{12}=\frac{\left\vert
\bigstar\left\{  1,2\right\}  \right\vert }{\left\vert \left\{  1,2\right\}
\right\vert }$ and $a_{34}=-\frac{\left\vert \bigstar\left\{  3,4\right\}
\right\vert }{\left\vert \left\{  3,4\right\}  \right\vert },$ we see that
$a_{12}<0$ and $a_{34}<0$ if and only if $\mathcal{T}^{\prime}$ is regular
after the flip and not regular before the flip. Since all areas $A_{ijk}$ are
positive, $a_{12}<0$ if and only if $\Phi<0$ and hence the result is proven.
\end{proof}

Note that in the proof we have shown that $\Phi<0$ if and only if
$\mathcal{T}$ is not regular and $\mathcal{T}^{\prime}$ is regular.

In order to get the global statement, one needs to know that a regular
triangulation can be found using flips. This is not true in general (see
\cite{ES}). However, we investigated some conditions when a flip algorithm
does work in Section \ref{flip algorithms}.

As a corollary of Rippa's theorem, we get an entropy quantity that increases
under the action of flipping to make a hinge regular.

\begin{corollary}
\label{eigenvalue cor}Consider the entropy defined by
\[
\Lambda=\inf\left\{  E\left(  f\right)  :\sum_{i\in\mathcal{T}_{0}}f_{i}%
^{2}=1\text{ and }\sum_{i\in\mathcal{T}_{0}}f_{i}=0\right\}  .
\]
Then $\Lambda$ decreases when an edge is flipped to make the hinge regular.
\end{corollary}

\begin{proof}
Let $\Lambda^{\prime}$ denote the entropy after the flip and let $f_{0}$ be
the $f$ which realize $\Lambda$ (since $f$ is in a compact set, there must be
an actual $f$ which minimizes $E\left(  f\right)  $). Then
\[
\Lambda^{\prime}=\inf_{f}E_{\mathcal{T}^{\prime}}\left(  f\right)  \leq
E_{\mathcal{T}^{\prime}}\left(  f_{0}\right)  \leq E_{\mathcal{T}}\left(
f_{0}\right)  =\Lambda.
\]

\end{proof}

Note that $\Lambda$ can be considered an eigenvalue of a particular operator
closely related to $\triangle.$ We remark that Corollary \ref{eigenvalue cor}
is similar in spirit to what is proven by G. Perelman at the beginning of his
paper \cite{Per}, where he shows that a slightly more complicated entropy,
\[
\inf\left\{  \int\left(  Rf^{2}+4\left\vert \nabla f\right\vert ^{2}\right)
dV:\int f^{2}dV=1\right\}  ,
\]
where $R$ is the scalar curvature, increases under Ricci flow.

Note that in $n$ dimensions, the regularity condition corresponds to
$\left\vert \bigstar\sigma^{n-1}\right\vert >0$ while good Dirichlet energy
corresponds to $\left\vert \bigstar\sigma^{1}\right\vert >0.$ Hence the
correspondence between regular triangulations and the Dirichlet energy only
occurs in dimension $2$ because $1=2-1,$ which is why the theorem is only
described for dimension $2.$ Although we do not pursue it here, this may
indicate that the Laplacian should instead be defined on functions on vertices
of the dual complex, $f:\bigstar\mathcal{T}_{n}\rightarrow\mathbb{R}$, in
which case the Laplacian would be
\[
\left(  \triangle f\right)  _{\bigstar\sigma_{0}^{n}}=\frac{1}{\left\vert
\sigma_{0}^{n}\right\vert }\sum_{\sigma^{n}\in\mathcal{T}_{n}}\frac{\left\vert
\sigma^{n}\cap\sigma_{0}^{n}\right\vert }{\left\vert \bigstar\left(
\sigma^{n}\cap\sigma_{0}^{n}\right)  \right\vert }\left(  f_{\bigstar
\sigma^{n}}-f_{\bigstar\sigma_{0}^{n}}\right)
\]
where the sum is over all $n$-simplices. In this case, positivity of the
coefficients corresponds to being regular.

\subsection{Laplace and heat equations\label{laplace and heat equations}}

Given a Laplace operator, we can now consider the standard elliptic and
parabolic equations, namely the Laplace equation
\begin{equation}
\triangle u=0 \label{laplace equation}%
\end{equation}
and the heat equation%
\begin{equation}
\frac{du}{dt}=\triangle u, \label{heat equation}%
\end{equation}
where the heat equation is an ordinary differential equation since $\triangle$
is a difference operator. A solution $u$ to the Laplace equation is called a
\emph{harmonic function}.

In order to study these equations, it will sometimes be easier to consider
$\triangle u=0$ as a matrix equation. We think of $u:\mathcal{T}%
_{0}\rightarrow\mathbb{R}$ as a vector and $\triangle$ corresponds to a matrix
$L$ whose off-diagonal pieces are
\[
L_{ij}=\frac{\left\vert \bigstar\left\{  i,j\right\}  \right\vert }{\left\vert
\left\{  i,j\right\}  \right\vert }%
\]
and whose diagonal pieces are%
\[
L_{ii}=-\sum_{j:\left\{  i,j\right\}  \in\mathcal{T}_{1}}\frac{\left\vert
\bigstar\left\{  i,j\right\}  \right\vert }{\left\vert \left\{  i,j\right\}
\right\vert }.
\]
Then one can write the Laplace equation as
\[
Lu=0.
\]
Note that if we wish to consider Poisson's equation
\begin{equation}
\triangle u=f \label{poisson equation}%
\end{equation}
then this is equivalent to
\[
Lu=fV
\]
where $\left(  fV\right)  _{i}=f_{i}\left\vert \bigstar\left\{  i\right\}
\right\vert .$ It is clear that $L$ has the constant functions $f_{i}=a$ (or
the vector $\left(  a,a,\ldots,a\right)  $) in the nullspace. If $\left\vert
\bigstar\left\{  i,j\right\}  \right\vert >0$ then we find the following.

\begin{theorem}
\label{laplace neg semidef 1}If $\left\vert \bigstar\left\{  i,j\right\}
\right\vert >0$ for all edges $\left\{  i,j\right\}  $ then $L$ is negative
semidefinite with nullspace spanned by the constant vectors.
\end{theorem}

\begin{proof}
In this case we have an $N\times N$ matrix $L$ with diagonal entries negative
and off-diagonal entries positive and with $\sum_{j=1}^{N}L_{ij}=0.$ We
reiterate an argument from \cite{CR}. Let $\left(  v_{1},\ldots,v_{N}\right)
$ be an eigenvector corresponding to $\lambda\geq0.$ We may assume that
$v_{1}>0$ is the maximum of $v_{i}.$ We wish to show that $v_{i}=v_{j}$ for
all $i,j.$ Observe%
\[
\lambda v_{1}=\sum_{i=1}^{N}L_{1i}v_{i}\leq\sum_{i=1}^{N}L_{1i}v_{1}=0.
\]
Equality holds if and only if $v_{i}=v_{1}$ for all $i.$
\end{proof}

\begin{corollary}
If $\left\vert \bigstar\left\{  i,j\right\}  \right\vert >0$ for all edges
$\left\{  i,j\right\}  $ then Poisson's equation has a solution for any $f$
such that
\[
\sum_{i\in\mathcal{T}_{0}}f_{i}V_{i}=0.
\]

\end{corollary}

This is the analogue of the smooth result that $\triangle u=f$ has a solution
if $\int_{M}fdV=0.$ One may also consider boundary conditions such as
Dirichlet and Neumann conditions. These cases for Delaunay triangulations in
two dimensions were studied by Bobenko and Springborn \cite{BS}.

The condition $\left\vert \bigstar\left\{  i,j\right\}  \right\vert >0$ is
obviously very important for the proof of Theorem \ref{laplace neg semidef 1}.
In two dimensions, this condition is equivalent to being regular by Corollary
\ref{edge regular duality}. It is not always necessary to assume $\left\vert
\bigstar\left\{  i,j\right\}  \right\vert >0$, as seen in the following
special cases.

Recall that in two dimensions, if a duality triangulation is edge-positive,
then the flip algorithm finds a regular triangulation (Theorem
\ref{2D edge flip algo}). For a similar set of two-dimensional triangulations,
the Laplacian is negative semidefinite.

\begin{theorem}
\label{laplace neg semidef 2}For any triangulation such that $d_{ij}>0$ for
all $\left(  i,j\right)  \in\mathcal{T}_{1}^{+},$ the Laplacian matrix $L$ is
negative semidefinite with nullspace spanned by the constant vectors.
\end{theorem}

We begin with a series of claims and an important lemma before beginning the
proof. We shall prove this by a sequence of claims. For all of the claims it
is assumed that the weights $d_{ij}$ are all positive. We shall use
$h_{ij}=d_{\pm}\left[  C\left(  \left\{  1,2,3\right\}  \right)  ,C\left(
\left\{  i,j\right\}  \right)  \right]  $ and $\gamma_{i}$ is the angle at
vertex $i.$ Consider only the $3\times3$ matrix $M$ corresponding to $\left\{
1,2,3\right\}  $ with entries $M_{ij}=h_{ij}/\ell_{ij}$ if $i\neq j$ and
$M_{ii}=-\sum_{j\neq i}M_{ij}.$

\begin{claim}
If $h_{ij}<0$ then $\gamma_{i}<\frac{\pi}{2}$ and $\gamma_{j}<\frac{\pi}{2}.$
\end{claim}

\begin{proof}
Let $k$ be the third vertex so that $\left\{  i,j,k\right\}  =\left\{
1,2,3\right\}  .$ We know that
\[
h_{ij}=\frac{d_{ik}-d_{ij}\cos\gamma_{i}}{\sin\gamma_{i}}%
\]
by formula \ref{center distance}). If $h_{ij}<0$ then $0<d_{ik}<d_{ij}%
\cos\gamma_{i}.$ Hence $\cos\gamma_{i}>0$ and $\gamma_{i}<\pi/2.$ We can also
express $h_{ij}$ as
\[
h_{ij}=\frac{d_{jk}-d_{ji}\cos\gamma_{j}}{\sin\gamma_{j}}%
\]
and follow the same logic.
\end{proof}

Thus only one $M_{ij}$ may be negative. Suppose it is $M_{12}.$

\begin{claim}
$M_{12}+M_{13}=\frac{\ell_{23}\left(  d_{12}\cos\gamma_{2}+d_{13}\cos
\gamma_{3}\right)  }{2A_{123}}.$
\end{claim}

\begin{proof}
We calculate%
\begin{align*}
M_{12}+M_{13}  &  =\frac{d_{23}-d_{21}\cos\gamma_{2}}{\ell_{12}\sin\gamma_{2}%
}+\frac{d_{32}-d_{31}\cos\gamma_{3}}{\ell_{13}\sin\gamma_{3}}\\
&  =\frac{\ell_{23}\left(  \ell_{23}-d_{21}\cos\gamma_{2}-d_{31}\cos\gamma
_{3}\right)  }{2A_{123}}%
\end{align*}
and finally we use that $\ell_{23}=\ell_{12}\cos\gamma_{2}+\ell_{13}\cos
\gamma_{2}.$
\end{proof}

\begin{claim}
$d_{12}\cos\gamma_{2}+d_{13}\cos\gamma_{3}>0.$
\end{claim}

\begin{proof}
If both $\gamma_{2}$ and $\gamma_{3}$ are less than or equal to $\pi/2$ then
this is clear (since both may not be equal to $\pi/2$). Since $M_{12}<0,$ and
hence $h_{12}<0,$ we can only have $\gamma_{3}>\pi/2.$ Since $h_{12}<0$ and
$h_{13}>0$ we have that
\[
\frac{d_{13}}{d_{12}}<\cos\gamma_{1}<\frac{d_{12}}{d_{13}}%
\]
so $d_{12}>d_{13}.$ Furthermore, since $\gamma_{1}+\gamma_{2}<\pi$ we have
that%
\[
0<-\cos\gamma_{3}=\cos\left(  \gamma_{1}+\gamma_{2}\right)  <\cos\gamma_{2}%
\]
so
\[
-d_{13}\cos\gamma_{3}<d_{12}\cos\gamma_{2}.
\]

\end{proof}

\begin{lemma}
\label{diagonal entries negative lemma}$M_{ii}<0.$
\end{lemma}

\begin{proof}
By the above argument, we know that $M_{11}=-M_{12}-M_{13}<0.$ Similar
arguments hold for the other coefficients.
\end{proof}

\begin{proof}
[Proof of Theorem \ref{laplace neg semidef 2}]It is sufficient to prove that
for any matrix $M_{ij},$ $1\leq i,j\leq3,$ is negative semidefinite. We know
that the vector $\left(  1,1,1\right)  $ is in the nullspace and we have
already shown in Lemma \ref{diagonal entries negative lemma} that the diagonal
entries are negative. Hence it is sufficient to show that the determinant of
the $2\times2$ submatrix $M_{ij},$ $1\leq i,j\leq2,$ is positive. We find that
the $2\times2$ determinant is equal to $M_{12}M_{13}+M_{12}M_{23}+M_{13}%
M_{23}.$ We compute the determinant to be equal to
\[
\frac{\left(  d_{13}h_{23}+d_{23}h_{13}\right)  \sin\gamma_{2}}{\ell_{12}%
\ell_{13}}%
\]
(to do this calculation, begin by writing the terms in the determinant using
formula (\ref{center distance}) choosing all of the denominators to contain
$\sin\gamma_{1}\sin\gamma_{2},$ then rearrange the terms using the facts that
$\gamma_{1}+\gamma_{2}+\gamma_{3}=\pi,$ $d_{ij}+d_{ji}=\ell_{ij},$ and
$\ell_{ij}=\ell_{ik}\cos\gamma_{i}+\ell_{jk}\cos\gamma_{k}$ several times and
finally recollecting $h_{23}$ and $h_{13}$ again using formula
(\ref{center distance})). Note that the determinant is symmetric in all
permutations in $1,2,3.$ We know by the claim above that two of the three
$h_{ij}$ must be positive, so choosing the two that are positive, we must have
that the determinant is positive. Hence the matrix is negative semidefinite.
\end{proof}

We consider $d_{ij}$ to be the length of a vector located at $i$ and in the
direction towards $j.$ Thus the condition $d_{ij}>0$ is like a positivity (or
Riemannian) condition for a metric (which measures the length of vectors) and
is thus a somewhat natural condition. The following is another result on
definiteness of the Laplacian with different assumptions.

\begin{theorem}
For a three-dimensional sphere packing triangulation, $L$ is negative
semidefinite with nullspace spanned by the constant vectors.
\end{theorem}

\begin{proof}
It is proven in \cite{G2} (see also \cite{Riv2}) that the matrix $A_{\left\{
1,2,3,4\right\}  }=\left(  \frac{\partial\alpha_{i}}{\partial r_{j}}\right)
_{1\leq i,j\leq4}$ is negative semidefinite with nullspace spanned by the
vector $\left(  r_{1},\ldots,r_{4}\right)  $. If we let $R_{\left\{
1,2,3,4\right\}  }$ be the diagonal matrix with $r_{i},$ $i=1,\ldots,4$ on the
diagonal, we see that
\[
L=\sum_{\sigma^{3}\in\mathcal{T}_{3}}\left(  R_{\sigma^{3}}A_{\sigma^{3}%
}R_{\sigma^{3}}\right)  _{E}.
\]
where $\left(  M_{\sigma^{3}}\right)  _{E}$ is the matrix extended by zeroes
to a $\left\vert \mathcal{T}_{0}\right\vert \times\left\vert \mathcal{T}%
_{0}\right\vert $ matrix so that the $\left(  M_{\sigma^{3}}\right)  _{E}$
acts on a vector $\left(  v_{1},\ldots,v_{\left\vert \mathcal{T}%
_{0}\right\vert }\right)  $ only on the coordinates corresponding to vertices
in $\sigma^{3}.$ Since $r_{i}>0$ for all $i\in\mathcal{T}_{0},$ it follows
that $L$ is negative semidefinite with nullspace spanned by $\left(
1,\ldots,1\right)  .$
\end{proof}

The importance of this result is it does not assume any positivity of the dual
area, which appears to be stronger than the assumption that $L$ is negative
definite. If $L$ is negative semi-definite with nullspace spanned by the
constant vector $\left(  1,\ldots,1\right)  $ then one can always solve the
Poisson equation for $f$ such that $\sum f_{i}A_{i}=0.$

The heat equation is an time-dependent, linear ordinary differential equation%
\[
\frac{du}{dt}=Lu
\]
whose short time existence is guaranteed by the existence theorem for ordinary
differential equations. One of the key properties of the heat equation is the
maximum principle, which says that the maximum decreases and the minimum
increases. This is true if $\left\vert \bigstar\left\{  i,j\right\}
\right\vert >0$.

\begin{theorem}
If $\left\vert \bigstar\left\{  i,j\right\}  \right\vert >0$ then for a
solution $u_{i}\left(  t\right)  $ of the heat equation, $u_{\max}\left(
t\right)  $ decreases and $u_{\min}\left(  t\right)  $ increases, where
$u_{\max}=\max\left\{  u_{i}:i\in\mathcal{T}_{0}\right\}  $ and $u_{\min}%
=\min\left\{  u_{i}:i\in\mathcal{T}_{0}\right\}  .$
\end{theorem}

\begin{proof}
The proof is standard and is simply that for any operator $Eu$ defined by
\[
\left(  Eu\right)  _{i}=\sum_{j\neq i}e_{ij}\left(  u_{j}-u_{i}\right)
\]
for some weights $e_{ij}>0,$ then $\left(  Eu\right)  _{i}<0$ if
$u_{i}=u_{\max}$ and $\left(  Eu\right)  _{i}>0$ if $u_{i}=u_{\min}.$
\end{proof}

Note that the maximum principle is not equivalent to $L$ being negative
semidefinite; it is a stronger condition and the proof uses that the
coefficients off the diagonal are positive. However, for certain functions
(geometric ones which are related to the coefficients of the Laplacian), it
may be possible to show that the maximum decreases and the minimum increases.
We call this a maximum principle for the function $f$ and we say that the
operator is \emph{parabolic-like} for the function $f.$ In \cite{G2} it is
proven that the sphere-packing case is parabolic-like for a curvature function
$K.$

\section{Toward discrete Riemannian manifolds\label{Riemannian}}

Much of this work arose out of an attempt to describe Riemannian manifolds
using piecewise Euclidean methods. In this final section, we try to describe
some of the work already done toward this end. There are two different
philosophies. One is to find analogues of the Riemannian setting. The idea is
to set up a framework on which variational-type arguments may be made
analogously to those in the smooth setting. The other is to actually
approximate smooth Riemannian geometry with discrete geometric structures. We
shall briefly consider both of these.

\subsection{Analogues of Riemannian geometry}

In this paper we gave a discrete operator on duality triangulations which, it
was argued, is an analogue of the Laplacian on a Riemannian manifold. This
gives rise also to a discrete heat equation, which is an ordinary differential
equation in this setting. It is not hard to imagine that similar arguments
give rise to Laplace-Beltrami operators on forms with the proper definition of
forms. A $k$-form can be defined to be an element of the dual space to the
vector space spanned by the $k$-dimensional simplices. There are also dual
$k$-forms which are elements of the dual space to the vector space spanned by
the duals of the $\left(  n-k\right)  $-dimensional simplices. Hirani
\cite{Hir} describes how to use duality information as we have described to
define the Hodge star operation, and thereby the Laplace-Beltrami operator on
these forms. One may then ask about an analogue of the Hodge theorem. This has
been studied somewhat by R. Hiptmair \cite{Hipt}. Study of the
Laplace-Beltrami operator on manifolds is also related to the study of the
Laplacian and harmonic analysis on metrized graphs and electrical networks
(see \cite{DS}, \cite{BF}, \cite{BR}).

Another important aspect of Riemannian geometry is the study of geodesics,
which we recall are locally length-minimizing curves. In the setting of
piecewise Euclidean manifolds, the geodesics are piecewise linear. One may
then ask many questions about geodesics, such as the number of closed
geodesics (see Pogorelov's work on quasi-geodesics on convex surfaces
\cite{Pog}) and the size of the cut locus to a basepoint, the locus of points
with two or more geodesics connecting it to the basepoint (see Miller-Pak
\cite{MP}). Many results on geodesics on piecewise Euclidean manifolds were
found by D. Stone \cite{Sto}, which lead him to some possible definitions of
curvature. The discrete geodesic problem for polytopes in $\mathbb{R}^{3}$ was
studied extensively in \cite{MMP}.

Much of modern Riemannian geometry is concerned with different notions of
curvature, such as sectional, Ricci, and scalar. In the piecewise Euclidean
setting, there are a number of definitions of curvatures, although it is still
somewhat an open question which ones are the proper ones for classification
purposes. Since the literature in this area is vast, we simply indicate some
of the principle works. D. Stone \cite{Sto} was successful in proving
analogues of the Cartan-Hadamard theorem (that negatively curved manifolds
have universal cover homeomorphic to $\mathbb{R}^{n}$) and Myer's theorem
(that positively curved manifolds are compact with a bound on the diameter) on
piecewise Euclidean manifolds using a quantity which he calls bounds on
sectional curvature. T. Regge introduced a notion of scalar curvature which is
described at each $\left(  n-2\right)  $-dimensional simplex as $2\pi$ minus
the sum of the dihedral angles at that simplex \cite{Reg}. This has been
widely studied as the so-called \textquotedblleft Regge
calculus\textquotedblright\ (see, for instance, \cite{Fro}, \cite{HW},
\cite{Ham}, \cite{ACM}). There are even some convergence results, which we
mention in the next section. Another potential curvature quantity in three
dimensions is described by Cooper and Rivin in \cite{CR}. They consider the
curvature at a vertex to be $4\pi$ minus the sum of the solid (or trihedral)
angles at the vertex. This curvature is certainly weaker than the curvature
introduced by Regge, but may be related to scalar curvature. It is possible
that the right curvature quantity will lead to a geometric flow which
simplifies geometry in a way similar to the way Ricci or Yamabe flow do in the
smooth category. This has been studies a bit in \cite{CL}, \cite{Luo},
\cite{G1}, \cite{G2}, and actually was the initial motivation for the
definitions of Laplacian described in this paper. Other applications of
discrete analogues of Riemannian geometry or geometric operators can be found
in \cite{BS}, \cite{ILTC}, \cite{MDSB}, \cite{Mer}, \cite{PP}, and
\cite{WGCTY}. In addition, techniques applying to metric spaces with sectional
curvature bounded in the sense of Alexandrov may apply (see \cite{BBI}).

\subsection{Approximating Riemannian geometry}

Another goal is to approximate Riemannian geometry by a discrete geometry such
as piecewise Euclidean triangulations. One would hope to be able to find
elements of Riemannian geometry such as Laplacian, Levi-Civita connection,
sectional curvature, scalar curvature, and so forth and not only have
analogous structures, but be able to show that as the triangulation gets finer
and finer, the discrete versions converge to the smooth versions. We mention
here some of the results which have been successful in this direction.

One of the most influential works is by Cheeger, M\"{u}ller, and Schrader, who
were able to relate discrete curvatures to Lipschitz-Killing curvatures
\cite{CMS}. The relevant discrete curvature is the sum certain angles and
volumes of hinges. In particular, the scalar curvature measure ($RdV$) is
concentrated on $\left(  n-2\right)  $-dimensional hinges in a triangulation,
and under a condition that the triangulation does not degenerate, they find
that the curvature quantity $2\pi$ minus the sum of the dihedral angles
multiplied by the volume of the $\left(  n-2\right)  $-dimensional hinge
converges to the scalar curvature measure. This version of scalar curvature is
also the one suggested by Regge \cite{Reg} and used extensively in the Regge
calculus. They prove convergence for each of the Lipschitz-Killing curvatures.
In addition, Barrett and Parker \cite{BP} proved a pointwise convergence of
piecewise-linear approximations of the Riemannian metric tensor and certain
types of tensor fields.

In regards to the Laplacian, some experimental work has been done by G. Xu
studying pointwise convergence of different discretized Laplace-Beltrami
operators to the smooth ones \cite{Xu1} \cite{Xu2}. Some of the
discretizations are the same or similar to those considered in this paper,
while some are not. On graphs (one-dimensional manifolds and generalizations),
it has been shown that the eigenvalues of the discrete Laplacians on metrized
graphs converge to the eigenvalues of the smooth Laplacian on a metrized graph
\cite{Fuj1} \cite{Fuj2} \cite{Fab}.

It was W. Thurston's idea to approximate the Riemann mapping between subsets
of $\mathbb{C}$ by mappings of circle packings. Such a discretization has been
shown to actually converge to the Riemann mapping \cite{RS}.

\begin{acknowledgement}
I\ would like to thank Herbert Edelsbrunner, Feng Luo, and Igor Pak for some
very helpful conversations related to parts of this paper.
\end{acknowledgement}

\end{document}